%% file: slant-cell-advection.tex
\documentclass[times]{elsarticle}
\pdfoutput=1
\usepackage{fullpage}
\usepackage{graphicx}
\usepackage{amsmath}
\usepackage{mathtools}
\usepackage{xcolor}
\usepackage{bm}
\usepackage{natbib}
\usepackage[hidelinks]{hyperref}
\usepackage{doi}
\usepackage[final,babel]{microtype}
\usepackage[utf8]{inputenc}
\usepackage[british]{babel}
\usepackage{csquotes}
\usepackage[T1]{fontenc}
\usepackage{siunitx}
\usepackage[font={small}]{caption}

\newcommand{\iu}{{i\mkern1mu}}

\newcommand{\revone}[1]{#1}
\newcommand{\revtwo}[1]{#1}
\newcommand{\revthree}[1]{#1}
\newcommand{\revother}[1]{#1}

\begin{document}

\begin{frontmatter}
\title{Multidimensional method-of-lines transport for atmospheric flows over steep terrain using arbitrary meshes}
\author[uor]{James Shaw\corref{cor}}
\ead{js102@zepler.net}
\author[uor]{Hilary Weller}
\author[uor]{John Methven}
\author[mo]{Terry Davies}

\cortext[cor]{Corresponding author}
\address[uor]{Department of Meteorology, University of Reading, Reading, United Kingdom}
\address[mo]{Met Office, Exeter, United Kingdom}

\input{abstract}

\end{frontmatter}

\input{intro}

\input{transportSchemes}

\section{Results}
\label{sec:results}

Three idealised numerical tests are performed to compare the accuracy of the cubicFit transport scheme with the multidimensional linear upwind scheme and with other transport schemes in the literature.  The first test transports a tracer horizontally on two-dimensional, highly-distorted terrain-following meshes, following Sch\"{a}r et al. \citep{schaer2002}.
The second is a new test case that modifies the velocity field and tracer position from the first test in order to challenge the stability and accuracy of the transport schemes near mountainous lower boundaries.
The third test evaluates the cubicFit scheme on \revone{hexagonal-icosahedral meshes} and cubed-sphere meshes with deformational flow on a spherical Earth, as specified by Lauritzen et al. \citep{lauritzen2012}.

We have implemented the cubicFit transport scheme and the numerical test cases using the OpenFOAM CFD library because it enables a like-for-like comparison between mesh types and transport schemes.  We provide source code archives for the OpenFOAM implementation of the cubicFit scheme \citep{atmosfoam}, the ASAM cut cell mesh generator \citep{asam_grid} and associated OpenFOAM converter \citep{gmv2openfoam}, and the hexagonal-icosahedral mesh generator \citep{geodesic-mesh}.  For the numerical test cases presented here we also supply the source code \citep{atmostests} and result data \citep{atmostests-data}.

\input{horizontalAdvection}

\input{mountainAdvection}
\input{deformationSphere}

\input{conclusion}

\section{Acknowledgements}
James Shaw acknowledges support from a PhD studentship funded jointly by NERC grant NE/K500860/1 and the University of Reading with CASE support from the Met Office.
We are grateful to the Leibniz Institute for Tropospheric Research for providing their cut cell mesh generator, \revother{to the three anonymous reviewers for their helpful questions,} and to Dr Shing Hing Man for his assistance with candidate polynomial generation.  We also thank Dr Tristan Pryer at the University of Reading for useful discussions about the cubicFit transport scheme.

\input{vonNeumann}
\input{spherical}

\bibliographystyle{elsarticle-num}
\bibliography{references}

\end{document}

%% file: abstract.tex
\begin{abstract}
Including terrain in atmospheric models gives rise to mesh distortions near the lower boundary that can degrade accuracy and challenge the stability of transport schemes.
Multidimensional transport schemes avoid splitting errors on distorted, \revone{arbitrary meshes}, and method-of-lines schemes have low computational cost because they perform reconstructions at fixed points.

This paper presents a multidimensional method-of-lines finite volume transport scheme, ``cubicFit'', which is designed to be numerically stable on arbitrary meshes.
Constraints derived from a von Neumann stability analysis are imposed during model initialisation to obtain stability and improve accuracy in distorted regions of the mesh, and near steeply-sloping lower boundaries.
Reconstruction calculations depend upon the mesh only, needing just one vector multiply per face per \revtwo{time-stage} irrespective of the velocity field.

The cubicFit scheme is evaluated using three, idealised numerical tests.  The first is a variant of a standard horizontal transport test on severely distorted terrain-following meshes.
The second is a new test case that assesses accuracy near the ground by transporting a tracer at the lower boundary over steep terrain on terrain-following meshes, cut-cell meshes, and new, slanted-cell meshes that do not suffer from severe time-step constraints associated with cut cells.
The third, standard test deforms a tracer in a vortical flow on \revone{hexagonal-icosahedral meshes} and cubed-sphere meshes.
In all tests, cubicFit is stable and largely insensitive to mesh distortions, and cubicFit results are more accurate than those obtained using a multidimensional linear upwind transport scheme.  The cubicFit scheme is second-order convergent regardless of mesh distortions.
\end{abstract}

%% file: intro.tex
\section{Introduction}

Numerical simulations of atmospheric flows solve equations of motion that result in the transport of momentum, temperature, water species and trace gases.  The numerical representation of Earth's terrain complicates the transport problem because the mesh is necessarily distorted \revtwo{next to} the lower boundary.
As new atmospheric models use increasingly fine mesh spacing, meshes are able to resolve steep, small-scale slopes.  Numerical schemes in operational weather forecast models can perform poorly over large mountain ranges, exhibiting small-scale numerical noise in momentum \citep{walko-avissar2008b}, temperature, humidity \citep{schaer2002} and potential vorticity fields \citep{hoinka-zaengl2004}, or even violating the Courant--Friedrich--Lewy stability constraint resulting in so-called `grid-point storms' \citep{webster2003}.
A transport scheme is desired that yields stable and accurate solutions, particularly near the surface where the weather affects us directly.
We present a new transport scheme which is numerically stable on arbitrary meshes and which is generally insensitive to mesh distortions created by steep slopes.  It has a low computational cost since most calculations are not repeated every time-step because they depend upon the mesh geometry only.

There are two main methods for representing terrain in atmospheric models: terrain-following layers and cut cells.
Both methods modify regular meshes to produce distorted meshes with cells that are aligned in columns.  Most operational models use terrain-following layers in which horizontal mesh surfaces are moved upwards to accommodate the terrain.  A vertical decay function is chosen so that mesh surfaces slope less steeply with increasing height.
The most straightforward is the linear decay function used by the basic terrain-following transform \citep{galchen-somerville1975} (also called the sigma coordinate), but many atmospheric models suffer from large numerical errors on such meshes \citep{schaer2002,klemp2011,eckermann2014}.
To reduce such errors, more complex decay functions have been developed so that mesh surfaces are smoother \citep{simmons-burridge1981,schaer2002,leuenberger2010,klemp2011}.

An alternative to terrain-following layers is the cut cell method.  Cut cell meshes are constructed by `cutting' a regular mesh with a piecewise-linear representation of the terrain.  New vertices are created where the terrain intersects mesh edges, and cell volumes that lie beneath the ground are removed.  Cut cell meshes can have arbitrarily small cells that impose severe time-step constraints on explicit transport schemes.  Several techniques have been developed to alleviate this problem, known as the `small-cell problem': small cells can be merged with adjacent cells \citep{yamazaki2016}, cell volumes can be artificially increased \citep{steppeler2002}, or an implicit scheme can be used near the ground with an explicit scheme used aloft \citep{jebens2011}.

Another method for avoiding the small-cell problem was proposed by Shaw and Weller \citep{shaw-weller2016} in which cell vertices are moved vertically so that they are positioned at the terrain surface.  We refer to this alternative method as the slanted cell method in order to distinguish it from the traditional cut cell method.  Slanted cell meshes do not suffer from arbitrarily small cells because the horizontal cell dimensions are not modified by the presence of terrain.

Smoothed terrain-following layers, cut cells and slanted cell methods all reduce the amount of mesh distortion but any mesh that represents sloping terrain must necessarily be distorted, at least near the ground.
Even when distortions are minimal, transport across mesh surfaces tends to be more common near steep slopes, and this misalignment between the flow and mesh surfaces increases numerical errors \citep{leonard1993,schaer2002,shaw-weller2016}.
A huge variety of transport schemes have been developed for atmospheric models, but few are able to account for distortions associated with steep terrain because they treat horizontal and vertical transport separately \citep{kent2014}, resulting in numerical errors called `splitting errors'.
Such errors can be reduced by explicitly accounting for transverse fluxes when combining fluxes \citep{leonard1996}, but splitting errors are still apparent in flows over steep terrain where meshes are highly distorted and metric terms in a terrain-following coordinate transform are large \citep{weller2017}.

Transport schemes are often classified as dimensionally-split or multidimensional.
Dimensionally-split schemes such as \citep{lin-rood1996,katta2015} calculate transport in each dimension separately before the flux contributions are combined.  Such schemes are computationally efficient and allow existing one-dimensional high-order methods to be used.
\revtwo{
Dimensionally-split schemes have only been used with quadrilateral meshes where dimensions are inherently separable.  Special treatment is required at the corners of cubed-sphere panels where local coordinates differ \citep{putman-lin2007,katta2015}.
For similar reasons, dimensionally-split schemes have only been used with terrain-following coordinate transforms and not cut cells.
}
Perhaps confusingly, dimensionally-split schemes are sometimes called multidimensional, too, because they use one-dimensional techniques for multidimensional transport.

Unlike dimensionally-split schemes, multidimensional schemes consider transport in two or three dimensions together.
There are several subclasses of multidimensional schemes that include
semi-Lagrangian finite volume schemes (also called conservative mesh remapping),
swept-area schemes (also called flux-form semi-Lagrangian, incremental remapping, or forward-in-time),
and method-of-lines schemes (also called Eulerian schemes).
Two-dimensional semi-Lagrangian finite volume schemes such as \citep{iske-kaeser2004,lauritzen2010} integrate over departure cells that are found by tracing backward the trajectories of cell vertices.  These schemes are conservative because departure cells are constructed so that there are no overlaps or gaps, which requires that cell areas are simply-connected domains \citep{lauritzen2011book}.
\revtwo{SLICE-3D is a three-dimensional semi-Lagrangian finite volume scheme for latitude-longitude meshes that applies separate conservative remappings in each dimension \citep{zerroukat-allen2012}.}
Swept area schemes such as \citep{lashley2002,skamarock-menchaca2010,lauritzen2011,thuburn2014} calculate the flux through a cell face by integrating over the upstream area that is swept out over one time-step.  Such schemes differ in their choice of area approximation, sub-grid reconstruction, and spatial integration method.
Because swept area schemes integrate over the reconstructed field, they typically require a matrix-vector multiply per face \revother{per time-stage} \citep{thuburn2014,skamarock-menchaca2010}.
Method-of-lines schemes such as \citep{weller2009,skamarock-gassmann2011} use a spatial discretisation to reduce the transport PDE to an ODE that is typically solved using a multi-stage time-stepping method.
\revother{A method-of-lines scheme using a spectral element reconstruction was recently developed to achieve accurate solutions near the surface of cut cell meshes \citep{steppeler-klemp2017}.}
Unlike semi-Lagrangian finite volume schemes, swept-area and method-of-lines schemes \revtwo{achieve conservation for small-scale rotational flows.
Such flows can twist the departure domain to such an extent that the domain intersects itself \citep{lauritzen2011}.  In two dimensions, a self-intersecting departure domain has a bowtie or hourglass shape}.
There are many more types of atmospheric transport schemes, but all can be classified according to their treatment of the three spatial dimensions.  A more comprehensive overview is presented by Lauritzen et al. \cite{lauritzen2014}.

For transport schemes that are ordinarily classified as `multidimensional', a further distinction ought to made between horizontally-multidimensional and three-dimensional schemes.
Most multidimensional schemes are only horizontally-multidimensional because, while the two horizontal dimensions are considered together, horizontal and vertical transport are still treated separately.
This separate treatment becomes less justifiable as atmospheric models are using increasingly fine horizontal mesh spacings that resolve small-scale steep slopes, resulting in greater mesh distortion and possible splitting errors \citep{kent2014}.
Three-dimensional schemes avoid any splitting errors over steep slopes, but only a few conservative three-dimensional schemes have been used in atmospheric models.
The multi-moment constrained finite volume scheme \citep{ii-xiao2009} is a three-dimensional scheme that has been used to simulate nonhydrostatic flows over orography with terrain-following coordinates on a $x$--$z$ plane \citep{li2013}.  Simulations of subcritical flow around a cylinder have also been performed on a three-dimensional hexahedral-prismatic hybrid mesh \citep{xie-xiao2016}.
The Multidimensional Positive Definite Advection Transport Algorithm (MPDATA) is another three-dimensional scheme that is suitable for arbitrary meshes.
It has been used on triangular unstructured meshes to simulate two-dimensional nonhydrostatic flows over orography \citep{smolarkiewicz-szmelter2011}, and in three-dimensional transport tests \citep{smolarkiewicz-szmelter2005}. 
\revone{Most recently, MPDATA has been extended to enable semi-implicit integrations of the compressible Euler equations on arbitrary meshes \citep{kuehnlein-smolarkiewicz2017}}.
The three-dimensional \revother{method-of-lines} scheme developed by Weller and Shahrokhi \citep{weller-shahrokhi2014} has been used in two-dimensional flows over orography on Cartesian $x$--$z$ planes with distorted meshes \citep{shaw-weller2016,weller2017}.
\revother{This finite volume scheme uses a moving least-squares reconstruction that makes it suitable for arbitrary meshes.
This least-squares approach has been applied previously to shallow water flows \citep{cuetofelgueroso2006}, aeronautic \citep{cuetofelgueroso2007} and porous media \citep{white2017} simulations.}

In this paper, we present a new multidimensional method-of-lines scheme, `cubicFit', that improves the stability of the Weller and Shahrokhi scheme \citep{weller-shahrokhi2014} and avoids all splitting errors.  To reconstruct values at cell faces, the scheme fits \revtwo{a multidimensional cubic polynomial over an upwind-biased stencil} using a least-squares approach.  The implementation uses constraints derived from a von Neumann stability analysis to select appropriate polynomial fits for stencils in highly-distorted mesh regions.  Almost all of the least-squares procedure depends upon the mesh geometry only and reconstruction weights can be precomputed without knowledge of the velocity field or tracer field.
Hence, \revthree{the computational cost of the cubicFit scheme is lower than most swept-area schemes, being more comparable to dimensionally-split schemes}, requiring only $n$ multiplies per cell face per \revtwo{time-stage} where $n$ is the size of the stencil.  \revtwo{Based on numerical experiments, the scheme is found to be conditionally stable up to maximum Courant numbers of about \num{1.3} to \num{3.3}.}

The remainder of this paper is organised as follows.
Section~\ref{sec:transport} starts by discretising the transport equation using a method-of-lines approach before describing the cubicFit transport scheme and a multidimensional linear upwind transport scheme.
Section~\ref{sec:results} evaluates the cubicFit scheme using three idealised numerical tests.
The first test follows Sch\"ar et al. \citep{schaer2002}, transporting a tracer horizontally above steep mountains on two-dimensional, highly-distorted terrain-following meshes.
The second is a new test case designed to assess numerical accuracy next to a mountainous lower boundary.  In this test, a tracer placed at the ground is transported over steep slopes by a terrain-following velocity field on terrain-following, cut cell and slanted cell meshes.
The third is a standard test of deformational flow on a single-layer spherical Earth, specified by Lauritzen et al. \citep{lauritzen2012}, which we use to assess the cubicFit transport scheme on \revone{hexagonal-icosahedral meshes} and cubed-sphere meshes.
Concluding remarks are made in section~\ref{sec:conclusion}.

%% file: transportSchemes.tex
\section{Transport schemes for arbitrary meshes}
\label{sec:transport}
The transport of a dependent variable $\phi$ in a prescribed, non-divergent velocity field $\mathbf{u}$ is given by the equation		
\begin{align}		
	\frac{\partial \phi}{\partial t} + \nabla \cdot \left( \mathbf{u} \phi \right) = 0 \text{ .} \label{eqn:advection}		
\end{align}
The time derivative is discretised using a two-stage, second-order Heun method,
\begin{subequations}
\begin{align}
	\phi^\star &= \phi^{(n)} + \Delta t \: g(\phi^{(n)}) \\
	\phi^{(n+1)} &= \phi^{(n)} + \frac{\Delta t}{2} \left[ g(\phi^{(n)}) + g(\phi^{\star}) \right]
\end{align} \label{eqn:heun}
\end{subequations}
where \(g(\phi^{(n)}) = - \nabla \cdot (\mathbf{u} \phi^{(n)})\) at time level \(n\).  The same time-stepping method is used for both the cubicFit scheme and the multidimensional linear upwind scheme.
Although the Heun method is unstable for a linear oscillator \citep{durran2013} and for solving the transport equation using centred, linear differencing, it is stable when it is used for transport schemes with \revtwo{sufficient upwinding}.

Using the finite volume method, the velocity field is prescribed at face centroids and the dependent variable is stored at cell centroids.  The divergence term in equation~\eqref{eqn:advection} is discretised using Gauss's theorem:
\begin{align}
	\nabla \cdot \left( \mathbf{u} \phi \right) \approx \frac{1}{\mathcal{V}_c} \sum_{f \in\:c} \mathbf{u}_f \cdot \mathbf{S}_f \phi_F \label{eqn:gauss-div}
\end{align}
where subscript $f$ denotes a value stored at a face and subscript $F$ denotes a value approximated at a face from surrounding values.  $\mathcal{V}_c$ is the cell volume, $\mathbf{u}_f$ is a velocity vector prescribed at a face, ${\mathbf{S}_f}$ is the surface area vector with a direction outward normal to the face and a magnitude equal to the face area, $\phi_F$ is an approximation of the dependent variable at the face, and $\sum_{f \in\:c}$ denotes a summation over all faces $f$ bordering cell $c$.
Note that equation~\eqref{eqn:gauss-div} is a second-order approximation of the divergence term which limits the cubicFit transport scheme to second-order numerical convergence.

This discretisation is applicable to arbitrary meshes.  A necessary condition for stability is given by the multidimensional Courant number,
\begin{align}
	\mathrm{Co}_c = \frac{\Delta t}{2 \mathcal{V}_c} \sum_{f \in\: c} \lvert \mathbf{u} \cdot \mathbf{S}_f \rvert \label{eqn:co}
\end{align}
such that $\mathrm{Co}_c \leq 1$ for all cells $c$ in the domain.  Hence, stability is constrained by the maximum Courant number of any cell in the domain.

The accurate approximation of the dependent variable at the face, $\phi_F$, is key to the overall accuracy of the transport scheme. The cubicFit scheme and multidimensional linear upwind scheme differ in their approximations, and these approximation \revone{methods are described} next.

\input{cubicFit}

\subsection{Multidimensional linear upwind transport scheme}
The multidimensional linear upwind scheme, called ``linearUpwind'' hereafter, is documented here since it provides a baseline accuracy for the experiments in section~\ref{sec:results}.  The approximation of $\phi_F$ is calculated using a gradient reconstruction,
\begin{align}
	\phi_F &= \phi_u + \nabla_c\: \phi \cdot \left(\mathbf{x}_f - \mathbf{x}_c \right)
\end{align} 
where $\phi_u$ is the upwind value of $\phi$, and $\mathbf{x}_f$ and $\mathbf{x}_c$ are the position vectors of the face centroid and cell centroid respectively.
The gradient $\nabla_c \:\phi$ is calculated using Gauss' theorem:
\begin{align}
	\nabla_c\: \phi = \frac{1}{\mathcal{V}_c} \sum_{f\in\:c} \widetilde{\phi}_F \mathbf{S}_f \label{eqn:linearUpwind-grad}
\end{align}
where $\widetilde{\phi}_F$ is linearly interpolated from the two neighbouring cells of face $f$.
\revone{The resulting stencil comprises all cells sharing a face with the upwind cell, including the upwind cell itself.  For a face in the interior of a two-dimensional rectangular mesh, the stencil for the linearUpwind scheme is a `$+$' shape with 5 cells.  On the same mesh, the stencil for the cubicFit scheme is more than twice the size with 12 cells.}
For cells adjacent to boundaries having zero gradient boundary conditions, the boundary value is set to be equal to the cell centre value before equation~\eqref{eqn:linearUpwind-grad} is evaluated.
This implementation of the multidimensional linear upwind scheme is included in the OpenFOAM software distribution \citep{openfoam}.

%% file: cubicFit.tex
\subsection{Cubic fit transport scheme}

\begin{figure}
	\centering
	\includegraphics{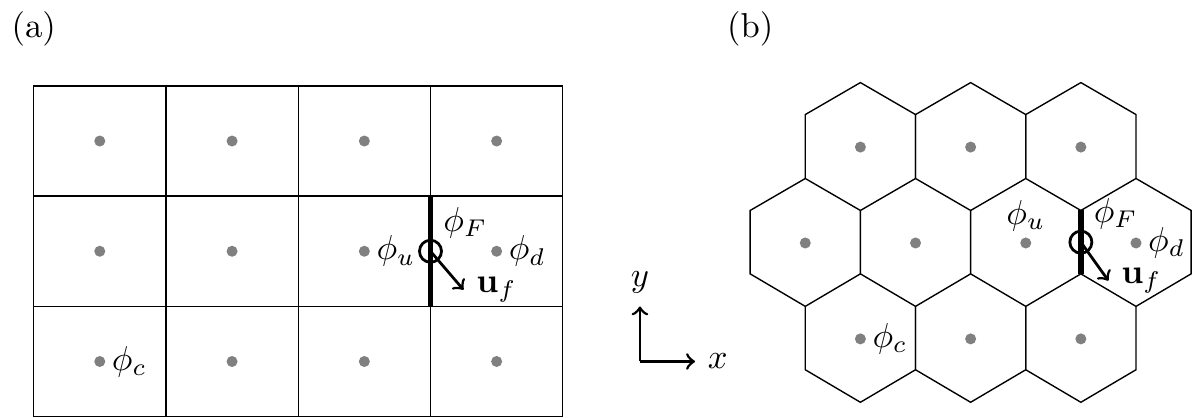}
	\caption{Upwind-biased stencils for faces far away from the boundaries of two-dimensional (a) rectangular and (b) hexagon meshes.  The stencil is used to fit a multidimensional polynomial to cell centre values, $\phi_c$, marked by grey circles, in order to approximate the value $\phi_F$ at the face centroid marked by an open circle.  $\phi_u$ and $\phi_d$ are the values at the centroids of the upwind and downwind cells neighbouring the target face, drawn with a heavy line.  The velocity vector $\mathbf{u}_f$ is prescribed at face $f$ and determines the choice of stencil at each time-step.}
	\label{fig:interiorStencils}
\end{figure}

The cubicFit scheme approximates the value of the dependent variable at the face, $\phi_F$, using a least-squares fit over a stencil of surrounding known values.
To introduce the approximation method, we will consider how an approximate value is calculated for a face that is far away from the boundaries of a two-dimensional uniform rectangular mesh.
For any mesh, every interior face connects two adjacent cells.  The velocity direction at the face determines which of the two adjacent cells is the upwind cell.  Since the stencil is upwind-biased and asymmetric, two stencils must be constructed for every interior face, and the appropriate stencil is chosen depending on the velocity direction at each face for every time-step.

The upwind-biased stencil for a face $f$ is shown in figure~\ref{fig:interiorStencils}a.  The wind at the face, $\mathbf{u}_f$, is blowing from the upwind cell $c_u$ to the downwind cell $c_d$.
To obtain an approximate value at $f$, a polynomial least-squares fit is calculated using the stencil values.
The stencil has \num{4} points in $x$ and \num{3} points in $y$, leading to a natural choice of polynomial that is cubic in $x$ and quadratic in $y$,
\begin{align}
	\phi = a_1 + a_2 x + a_3 y + a_4 x^2 + a_5 xy + a_6 y^2 + a_7 x^3 + a_8 x^2 y + a_9 x y^2 \label{eqn:fullPoly} \text{ .}
\end{align}
A least-squares approach is needed because the system of equations is overconstrained, with \num{12} stencil values but only \num{9} polynomial terms.  The stencil geometry is expressed in a local coordinate system with the face centroid as the origin so that the approximated value $\phi_F$ is equal to the constant coefficient $a_1$.
\revtwo{The stencil is upwind-biased to improve numerical stability, and the multidimensional cubic polynomial is chosen to improve accuracy in the direction of flow \citep{leonard1993}.}

The remainder of this subsection generalises the approximation technique for arbitrary meshes and describes the methods for constructing stencils, performing a least-squares fit with a suitable polynomial, and ensuring numerical stability of the transport scheme.

\subsubsection{Stencil construction}
\label{sec:stencil}

For every interior face, two stencils are constructed, one for each of the possible upwind cells.
Stencils are not constructed for boundary faces because values of $\phi$ at boundaries are calculated from prescribed boundary conditions.
For a given interior face $f$ and upwind cell $c_u$, we find those faces that are connected to $c_u$ and `oppose' face $f$.  These are called the \textit{opposing faces}.
The opposing faces for face $f$ and upwind cell $c_u$ are determined as follows.
Defining $G$ to be the set of faces other than $f$ that border cell $c_u$, we calculate the `opposedness', $\mathrm{Opp}$, between faces $f$ and $g \in G$, defined as
\begin{align}
	\mathrm{Opp}(f, g) \equiv - \frac{\mathbf{S}_f \cdot \mathbf{S}_g}{|\mathbf{S}_f|^2} \label{eqn:opp}
\end{align}
where $\mathbf{S}_f$ and $\mathbf{S}_g$ are the surface area vectors pointing outward from cell $c_u$ for faces $f$ and $g$ respectively.
Using the fact that $\mathbf{a} \cdot \mathbf{b} = |\mathbf{a}|\:|\mathbf{b}| \cos(\theta)$ we can rewrite equation~\eqref{eqn:opp} as
\begin{align}
	\mathrm{Opp}(f, g) = - \frac{|\mathbf{S}_g|}{|\mathbf{S}_f|} \cos(\theta)
\end{align}
where $\theta$ is the angle between faces $f$ and $g$.  In this form, it can be seen that $\mathrm{Opp}$ is a measure of the relative area of $g$ and how closely it parallels face $f$.

The set of opposing faces, $\mathrm{OF}$, is a subset of $G$, comprising those faces with $\mathrm{Opp} \geq 0.5$, and the face with the maximum opposedness.  Expressed in set notation, this is
\begin{align}
	\mathrm{OF}(f,c_u) \equiv \{ g : \mathrm{Opp}(f, g) \geq 0.5 \} \cup \{ g : \max_{g\:\in\:G}(\mathrm{Opp}(f, g)) \} \text{ .}
\end{align}
On a rectangular mesh, there is always one opposing face $g$, and it is exactly parallel to the face $f$ such that $\mathrm{Opp}(f, g) = 1$.

Once the opposing faces have been determined, the set of internal and external cells must be found.  The \textit{internal cells} are those cells that are connected to the opposing faces.  Note that $c_u$ is always an internal cell.  The \textit{external cells} are those cells that share vertices with the internal cells.  Note that $c_d$ is always an external cell.  Finally, the \textit{stencil boundary faces} are boundary faces having Dirichlet boundary conditions\footnote{Boundary faces with Neumann boundary conditions would require extrapolated boundary values to be calculated.
This would create a feedback loop in which boundary values are extrapolated from interior values, then interior values are transported using stencils that include boundary values.  
We have not considered how such an extrapolation \revtwo{could be made} consistent with the multidimensional polynomial reconstruction.
Hence, boundary faces with Neumann boundary conditions are excluded from the set of stencil boundary faces.} that share a vertex with the internal cells.
Having found these three sets, the stencil is constructed to comprise all internal cells, external cells and stencil boundary faces.

\begin{figure}
	\centering
	\includegraphics{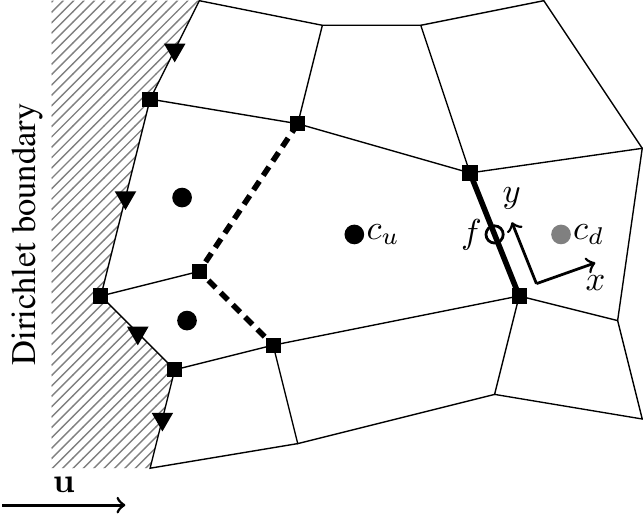}
	\caption{A fourteen-point, upwind-biased stencil for face $f$ connecting the pentagonal upwind cell, $c_u$, and the downwind cell $c_d$.  The dashed lines denote the two faces of cell $c_u$ that oppose $f$, and black circles mark the centroids of the internal cells that are connected to these two opposing faces.  The stencil is extended outwards by including cells that share vertices with the three internal cells, where black squares mark these vertices.  Four stencil boundary faces, marked by black triangles, are also included.
The local coordinate system $(x, y)$ has its origin at the centroid of face $f$, marked by an open circle, with $x$ normal to $f$ and $y$ perpendicular to $x$.}
	\label{fig:double-upwind-stencil}
\end{figure}

Figure~\ref{fig:double-upwind-stencil} illustrates a stencil construction for face $f$ connecting upwind cell $c_u$ and downwind cell $c_d$.  The two opposing faces are denoted by thick dashed lines and the centres of the three adjoining internal cells are marked by black circles.  The stencil is extended outwards by including the external cells that share vertices with the internal cells, where the vertices are marked by black squares.  A boundary at the far left has \revone{Dirichlet} boundary conditions, and so the four stencil boundary faces are also included in the stencil, where the boundary face centres are marked by black triangles.  The resultant stencil contains fourteen points.

\subsubsection{Least-squares fit}
To approximate the value of $\phi$ at a face $f$, a least-squares fit is calculated from a stencil of surrounding known values.  First, we will show how a polynomial least-squares fit is calculated for a face on a rectangular mesh.  Second, we will make modifications to the least-squares fit that are necessary for numerical stability.  

For faces that are far away from the boundaries of a rectangular mesh, we fit the multidimensional polynomial given by equation~\eqref{eqn:fullPoly} that has nine unknown coefficients, $\mathbf{a} = a_1 \ldots a_9$, using the twelve cell centre values from the upwind-biased stencil, $\bm{\phi} = \phi_1 \ldots \phi_{12}$.  This yields a matrix equation
\begin{align}
	\begin{bmatrix}
		1 & x_1 & y_1 & x_1^2 & x_1 y_1 & y_1^2 & x_1^3 & x_1^2 y_1 & x_1 y_1^2 \\
		1 & x_2 & y_2 & x_2^2 & x_2 y_2 & y_2^2 & x_2^3 & x_2^2 y_2 & x_2 y_2^2 \\
		\vdots & \vdots & \vdots & \vdots & \vdots & \vdots & \vdots & \vdots & \vdots \\
		1 & x_{12} & y_{12} & x_{12}^2 & x_{12} y_{12} & y_{12}^2 & x_{12}^3 & x_{12}^2 y_{12} & x_{12} y_{12}^2 \\
	\end{bmatrix}
	\begin{bmatrix}
		a_1 \\
		a_2 \\
		\vdots \\
		a_9
	\end{bmatrix}
	=
	\begin{bmatrix}
		\phi_1 \\
		\phi_2 \\
		\vdots \\
		\phi_{12}
	\end{bmatrix}
\end{align}
which can be written as
\begin{align}
	\mathbf{B} \mathbf{a} = \bm{\phi} \label{eqn:unweightedLeastSquares} \text{ .}
\end{align}
The rectangular matrix $\mathbf{B}$ has one row for each cell in the stencil and one column for each term in the polynomial.  $\mathbf{B}$ is called the \textit{stencil matrix}, and it is constructed using only the mesh geometry.
A local coordinate system is established in which $x$ is normal to the face $f$ and $y$ is perpendicular to $x$.
The coordinates $(x_i, y_i)$ give the position of the centroid of the $i$th cell in the stencil.
A two-dimensional stencil is also used for the tests on spherical meshes in section~\ref{sec:deformationSphere}.  In these tests, \revthree{cell centres are projected perpendicular to a tangent plane at the face centre.  Previous studies found that results were largely insensitive to the projection method \citep{skamarock-gassmann2011,lashley2002}.}

The unknown coefficients $\mathbf{a}$ are calculated using the pseudo-inverse, $\mathbf{B}^+$, found by singular value decomposition,
\begin{align}
	\mathbf{a} = \mathbf{B}^+ \bm{\phi} \text{ .}
\intertext{Recall that the approximate value $\phi_F$ is equal to the constant coefficient $a_1$, which is a weighted mean of $\bm{\phi}$,} 
	a_1 = \begin{bmatrix}
		b_{1,1}^+ \\
		b_{1,2}^+ \\
		\vdots \\
		b_{1,12}^+
	\end{bmatrix}
	\cdot
	\begin{bmatrix}
		\phi_1 \\
		\phi_2 \\
		\vdots \\
		\phi_{12}
	\end{bmatrix} \label{eqn:cubicFit:weighted-sum}
\end{align}
where the weights $b_{1,1}^+ \ldots b_{1,12}^+$ are the elements of the first row of $\mathbf{B}^+$.
Note that the majority of the least-squares fit procedure depends on the mesh geometry only.  An implementation may precompute the pseudo-inverse for each stencil during model initialisation, and only the first row needs to be stored.  Since each face has two possible stencils depending on the orientation of the velocity relative to the face, the implementation stores two sets of weights for each face.
Knowledge of the values of $\bm{\phi}$ is only required to calculate the weighted mean given by equation~\eqref{eqn:cubicFit:weighted-sum}, which is evaluated once per face per \revtwo{time-stage}.

In the least-squares fit presented above, all stencil values contributed equally to the polynomial fit.
It is necessary for numerical stability that the polynomial fits the cells connected to face $f$ more closely than other cells in the stencil, as shown by \citep{lashley2002,skamarock-menchaca2010}.
To achieve this, we allow each cell to make an unequal contribution to the least-squares fit.
We assign an integer \textit{multiplier} to each cell in the stencil, $\mathbf{m} = m_1 \ldots m_{12}$, and multiply equation~\eqref{eqn:unweightedLeastSquares} to obtain
\begin{align}
	\mathbf{\tilde{B}} \mathbf{a} = \mathbf{m} \cdot \bm{\phi}
\end{align}
where $\mathbf{\tilde{B}} = \mathbf{M} \mathbf{B}$ and $\mathbf{M} = \mathrm{diag}(\mathbf{m})$.  The constant coefficient $a_1$ is calculated from the pseudo-inverse, $\mathbf{\tilde{B}}^+$,
\begin{align}
	a_1 = \mathbf{\tilde{b}_1^+} \cdot \mathbf{m} \cdot \bm{\phi} \label{eqn:weightedPinv}
\end{align}
where $\mathbf{\tilde{b}_1^+} = \tilde{b}_{1,1}^+ \ldots \tilde{b}_{1,12}^+$ are the elements of the first row of $\mathbf{\tilde{B}}^+$.
Again, $a_1$ is a weighted mean of $\bm{\phi}$, where the weights are now $\mathbf{\tilde{b}_1^+} \cdot \mathbf{m}$.  Values for $\mathbf{m}$ are chosen so that the cells connected to face $f$ make a greater contribution to the least-squares fit, as discussed later in section~\ref{sec:stabilisation}.

For faces of a non-rectangular mesh, or faces that are near a boundary, the number of stencil points and number of polynomial terms may differ: a stencil will have one or more cells and, for two-dimensional meshes, its polynomial will have between one and nine terms.  Additionally, the polynomial cannot have more terms than its stencil has cells because this would lead to an underconstrained system of equations.  The procedure for choosing suitable polynomials is discussed next.

\subsubsection{Polynomial generation}
\label{sec:polyCandidates}
\begin{figure}
	\centering
	\includegraphics{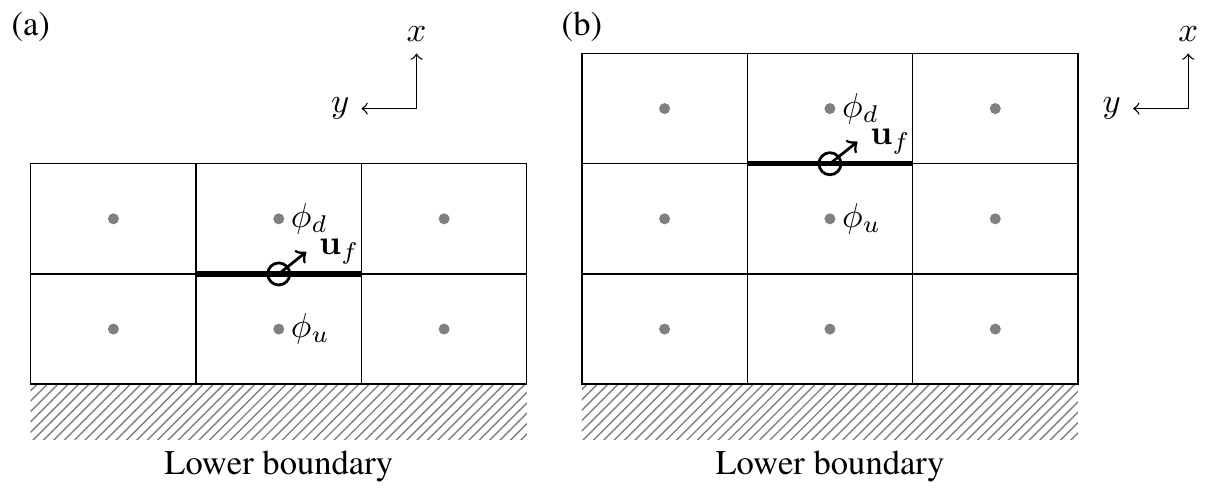}
	\caption{Upwind-biased stencils for faces near the lower boundary of a rectangular $x$--$z$ mesh, with (a) a $3 \times 2$ stencil for the face immediately adjacent to the lower boundary, and (b) a $3 \times 3$ stencil for the face immediately adjacent to the face in (a).  Each stencil belongs to the face marked by a thick line.  The local coordinate system is shown, having an $x$ direction normal to the face a $y$ direction tangent to the face.  For both stencils, attempting a least-squares fit using the nine-term polynomial in equation~\eqref{eqn:fullPoly} would result in an underconstrained problem.
	\revone{There is no normal flow at the lower boundary.}}
	\label{fig:boundaryStencils}
\end{figure}

The majority of faces on a uniform two-dimensional mesh have stencils with more than nine cells.  For example, a rectangular mesh has 12 points (figure~\ref{fig:interiorStencils}a), and a hexagonal mesh has 10 points (figure~\ref{fig:interiorStencils}b).
In both cases, constructing a system of equations using the nine-term polynomial in equation~\eqref{eqn:fullPoly} leads to an overconstrained problem that can be solved using least-squares.  However, this is not true for faces near boundaries: stencils that have fewer than nine cells (figure~\ref{fig:boundaryStencils}a) would result in an underconstrained problem, and stencils that have exactly nine cells may lack sufficient information to constrain high-order terms.  For example, the stencil in figure~\ref{fig:boundaryStencils}b lacks sufficient information to fit the $x^3$ term.  In such cases, it becomes necessary to perform a least-squares fit using a polynomial with fewer terms.

For every stencil, we find a set of \textit{candidate polynomials} that do not result in an underconstrained problem.
In two dimensions, a candidate polynomial has some combination of between one and nine terms from equation~\eqref{eqn:fullPoly}.  There are two additional constraints that a candidate polynomial must satisfy.

First, high-order terms may be included in a candidate polynomial only if the lower-order terms are also included.
More precisely, let
\revother{
\begin{align}
	M(x, y) = { x^i y^j : i,j \geq 0 \text{ and } i \leq 3 \text{ and } j \leq 2 \text{ and } i+j \leq 3}
\end{align}
be the set of all monomials of degree at most \num{3} in $x, y$.
A subset $S$ of $M(x,y)$ is ``dense'' if, whenever $x^a y^b$ is in $S$, then $x^i y^j$ is also in $S$ for all $0 \leq i \leq a$, $0 \leq j \leq b$.}
For example, the polynomial $\phi = a_1 + a_2 x + a_3 y + a_4 xy + a_5 x^2 + a_6 x^2 y$ is a dense subset of $M(x,y)$, but $\phi = a_1 + a_2 x + a_3 y + a_4 x^2 y$ is not because $x^2 y$ can be included only if $xy$ and $x^2$ are also included.
\revother{In total there are 26 dense subsets of the two-dimensional polynomial in equation~\eqref{eqn:fullPoly}.}

Second, a candidate polynomial must have a stencil matrix $\mathbf{B}$ that is full rank.  The matrix is considered full rank if its smallest singular value is greater than \num{1e-9}.
Using a polynomial with all nine terms and the stencil in figure~\ref{fig:boundaryStencils}b results in a rank-deficient matrix and so the nine-term polynomial is not a candidate polynomial.

The candidate polynomials are all the dense subsets of $M(x,y)$ \revother{that have a cardinality greater than one} with a stencil matrix that is full rank.  The final stage of the cubicFit transport scheme selects a candidate polynomial and ensures that the least-squares fit is numerically stable.

\subsubsection{Stabilisation procedure}
\label{sec:stabilisation}
So far, we have constructed a stencil and found a set of candidate polynomials.  Applying a least-squares fit to any of these candidate polynomials avoids creating an underconstrained problem.  The final stage of the transport scheme chooses a suitable candidate polynomial and appropriate multipliers $\mathbf{m}$ so that the fit is numerically stable.

The approximated value $\phi_F$ is equal to $a_1$ which is calculated from equation~\eqref{eqn:weightedPinv}.  The value of $a_1$ is a weighted mean of $\bm{\phi}$ where $\mathbf{w} = \mathbf{\tilde{b}_1^+} \cdot \mathbf{m}$ are the weights.
If the cell centre values $\bm{\phi}$ are assumed to approximate a smooth field then we expect $\phi_F$ to be close to the values of $\phi_u$ and $\phi_d$, and expect $\phi_F$ to be insensitive to small changes in $\bm{\phi}$.  When the weights $\mathbf{w}$ have large magnitude then this is no longer true: $\phi_F$ becomes sensitive to small changes in $\bm{\phi}$ which can result in large, \revtwo{numerically unstable} departures from the smooth field $\bm{\phi}$.

\revtwo{To avoid numerical instabilities, a simplified, one-dimensional von Neumann analysis was performed, presented in appendix A.  The analysis is used to impose three stability conditions on the weights $\mathbf{w}$},
\begin{subequations}
\label{eqn:stability}
\begin{align}
	0.5 \leq w_u \leq 1 \label{eqn:stabilityU} \\
	0 \leq w_d \leq 0.5 \label{eqn:stabilityD} \\
	w_u - w_d \geq \max_{p\:\in\:P}(|w_p|)
\end{align}
\end{subequations}
where $w_u$ and $w_d$ are the weights for the upwind and downwind cells respectively.  The \textit{peripheral points} $P$ are the cells in the stencil that are not the upwind or downwind cells, and $w_p$ is the weight for a given peripheral point $p$.
 The upwind, downwind and peripheral weights sum to one such that $w_u + w_d + \sum_{p \in P} w_p = 1$.

The stabilisation procedure comprises three steps.  In the first step, the set of candidate polynomials is sorted in preference order so that candidates with more terms are preferred over those with fewer terms.
If there are multiple candidates with the same number of terms, the minimum singular value of $\mathbf{B}$ is calculated for each candidate, and an ordering is imposed such that the candidate with the larger minimum singular value is preferred.  This ordering ensures that the preferred candidate is the highest-order polynomial with the most information content.\footnote{\revthree{Note that singular values are used for two purposes: first, to test if the matrix $\mathbf{B}$ is full-rank and, second, to impose an ordering on candidates.
We have used the minimum singular value, $\sigma_\mathrm{min}(\mathbf{B})$, for both purposes.  Alternatively, we could use the condition number, $\mathrm{cond}(\mathbf{B})$, which is the ratio of smallest to largest singular value.
Experiments revealed that only the candidate ordering was sensitive to the choice of $\sigma_\mathrm{min}$ or $\mathrm{cond}$.
The most suitable choices of singular value calculations could be explored in future.}}

In the second step, the most-preferred polynomial is taken from the list of candidates and the multipliers are assigned so that the upwind cell and downwind cell have multipliers $m_u = 2^{10}$ and $m_d = 2^{10}$ respectively, and all peripheral points have multipliers $m_p = 1$.  These multipliers are very similar to those used by \citep{lashley2002}, leading to a well-conditioned matrix $\mathbf{\tilde{B}}$ and a least-squares fit in which the polynomial passes almost exactly through the upwind and downwind cell centre values.

In the third step, we calculate the weights $\mathbf{w}$ and evaluate them against the stability \revother{conditions} given in equation~\eqref{eqn:stability}.  If any \revother{condition} is violated, the value of $m_d$ is halved and the \revother{conditions} are evaluated with the new weights.  This step is repeated until the weights satisfy the stability \revother{conditions}, or $m_d$ becomes smaller than one.  In practice, the \revother{conditions} are satisified when $m_d$ is either small (between 1 and 4) or equal to $2^{10}$.  The upwind multiplier $m_u$ is fixed at $2^{10}$ and the peripheral multipliers $m_p$ are fixed at \num{1}.  If the \revother{conditions} are still not satisfied, then we start again from the second step with the next polynomial in the candidate list. 

Finally, if no stable weights are found for any candidate polynomial, we revert to an upwind scheme such that $w_u = 1$ and all other weights are zero.
In our experiments we have not encountered any stencil for which this last resort is required.
\revone{Furthermore, our experiments show that the stabilisation procedure only modifies the least squares fit for stencils near boundaries and for stencils in distorted mesh regions.
For stencils in the interior of a uniform rectangular mesh, the least squares fit includes all terms in equation~\eqref{eqn:fullPoly} with $m_u = m_d = 2^{10}$.}

\begin{figure}
	\centering
	\includegraphics{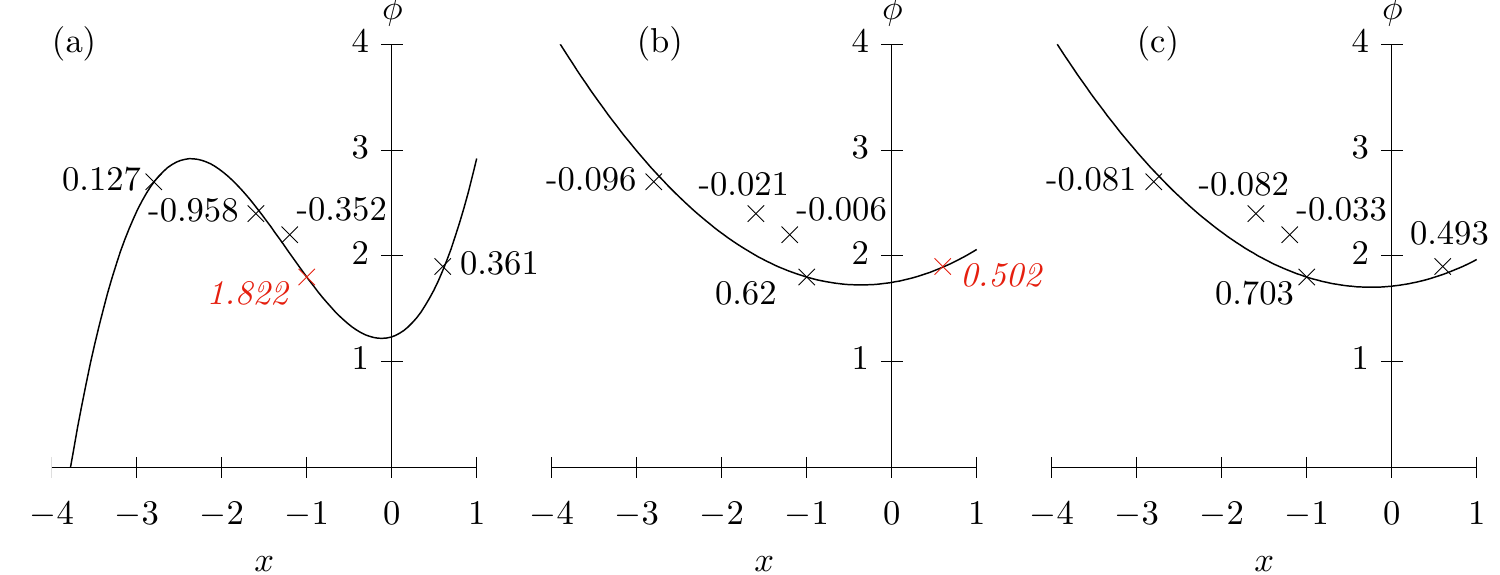}
	\caption{\revone{One-dimensional least-squares fits with} a stencil of five points using (a) a cubic polynomial with multipliers $m_u = 1024$, $m_d = 1024$ and $m_p = 1$, (b) a quadratic polynomial with the same multipliers, and (c) a quadratic polynomial with multipliers $m_u = 1024$, $m_d = 1$ and $m_p = 1$.  Notice that the curves in (a) and (b) fit almost exactly through the upwind and downwind points immediately adjacent to the $y$-axis, but in (c) the curve fits almost exactly only through the upwind point immediately to the left of the $y$-axis.  The point data are labelled with their respective weights.  Points that have failed one of the stability \revother{conditions} in equation~\eqref{eqn:stability} are marked in red with italicised labels.  The upwind point is located at $(-1, 1.8)$ and the downwind point at $(0.62, 1.9)$, and the peripheral points are at $(-2.8, 2.4)$, $(-1.6, 2.7)$ and $(-1.2, 2.2)$.  The stabilisation procedure (section~\ref{sec:stabilisation}) calculates weights using only $x$ positions, and values of $\phi$ are included here for illustration only.}
	\label{fig:oscillatory1D}
\end{figure}

To illustrate the stabilisation procedure, figure~\ref{fig:oscillatory1D}a presents a one-dimensional example of a cubic polynomial fitted through five points, with the weight at each point printed beside it.
The stabilisation procedure only uses the $x$ positions of these points and does not use the values of $\phi$ themselves.  The $\phi$ values are included here for illustration only.
Hence, for a given set of $x$ positions, the same set of weights are chosen irrespective of the $\phi$ values.

For a one-dimensional cubic polynomial fit, the list of candidate polynomials in preference order is
\begin{align}
	\phi &= a_1 + a_2 x + a_3 x^2 + a_4 x^3 \label{eqn:cubicCandidate} \text{ ,} \\
	\phi &= a_1 + a_2 x + a_3 x^2 \label{eqn:quadCandidate} \text{ ,} \\
	\phi &= a_1 + a_2 x \text{ ,} \\
	\phi &= a_1 \text{ .}
\end{align}
We begin with the cubic equation~\eqref{eqn:cubicCandidate}.  The multipliers are chosen so that the polynomial passes almost exactly through the upwind and downwind points that are immediately to the left and right of the $y$-axis respectively.
The \revother{stability condition} on the upwind point is violated because $w_u = 1.822 > 1$ (equation~\ref{eqn:stabilityU}).  Reducing the downwind multiplier does not help to satisfy the \revother{stability condition}, so we start again with the quadratic equation~\eqref{eqn:quadCandidate}, and the new fit is presented in figure~\ref{fig:oscillatory1D}b.
Again, the multipliers are chosen to force the polynomial through the upwind and downwind points, but this violates the \revother{stability condition} on the downwind point because $w_d = 0.502 > 0.5$ (equation~\ref{eqn:stabilityD}).  This time, however, stable weights are found by reducing \revtwo{$m_d$} to one (figure~\ref{fig:oscillatory1D}c) and these are the weights that will be used to approximate $\phi_F$, where the polynomial intercepts the $y$-axis.

\subsubsection{Future extension to three dimensions}
\revone{All the procedures used in the cubicFit scheme generalise to three dimensions.  The stencil construction procedure described in section~\ref{sec:stencil} creates a stencil with \num{12} cells for a face in the interior of a two-dimensional rectangular mesh.  In three dimensions, the same procedure creates a stencil with $3 \times 12 = 36$ cells.
A three-dimensional stencil has three times as many cells as its two-dimensional counterpart if the mesh has prismatic cells arranged in columns.  Hence, the computational cost during integration increases three-fold when moving from two dimensions to three dimensions.}

\revone{To extend the least squares fit to three dimensions, the two-dimensional polynomial in equation~\eqref{eqn:fullPoly} is replaced with its three-dimensional counterpart,
\begin{multline}
	\phi = a_1 + a_2 x + a_3 y + a_4 z + a_5 x^2 + a_6 xy + a_7 y^2 + a_8 xz + a_9 yz + a_{10} z^2 + \\ a_{11} x^3 + a_{12} x^2 y + a_{13} x y^2 + a_{14} x^2 z + a_{15} x z^2 + a_{16} y z^2 + a_{17} y^2 z + a_{18} xyz \text{ .} \label{eqn:fullPoly3D}
\end{multline}
The procedure for generating candidate polynomials described in section~\ref{sec:polyCandidates} results in 26 dense subsets in two dimensions and 842 dense subsets in three dimensions.  Note that the combinatorial explosion of dense subsets in three dimensions does not increase the computational cost during integration.}

\revone{The stabilisation procedure described in section~\ref{sec:stabilisation} requires further numerical experiments to verify that it is sufficient for three-dimensional flows and arbitrary polyhedral meshes.
An initial three-dimensional test with uniform flow and a uniform Cartesian mesh obtained a numerically stable result.
For stencils in the interior of the domain, the least squares fit includes all polynomial terms in equation~\eqref{eqn:fullPoly3D} with $m_u = m_d = 2^{10}$.
The stabilisation procedure does not modify the least squares fit for these stencils, but we have not explored the three-dimensional extension of cubicFit any further.}

%% file: horizontalAdvection.tex
\subsection{Horizontal transport over mountains}
\label{sec:horizontal}

A two-dimensional transport test was developed by Sch\"{a}r et al. \citep{schaer2002} to study the effect of terrain-following coordinate transformations on numerical accuracy.  In this standard test, a tracer is positioned aloft and transported horizontally over wave-shaped mountains.  The test challenges transport schemes because the tracer must cross mesh layers, which acts to reduce numerical accuracy \citep{schaer2002,shaw-weller2016}.
Here we use a more challenging variant of this test that has steeper mountains and highly-distorted terrain-following meshes.
Convergence results are compared using the linearUpwind and cubicFit transport schemes.

The domain is defined on a rectangular $x$--$z$ plane that is \SI{301}{\kilo\meter} wide and \SI{25}{\kilo\meter} high as measured between parallel boundary edges.
Boundary conditions are imposed on the tracer density $\phi$ such that $\phi = \SI{0}{\kilo\gram\per\meter\cubed}$ at the inlet boundary, and a zero normal gradient
$\partial \phi / \partial n = \SI{0}{\kilo\gram\per\meter\tothe{4}}$ is imposed at the outlet boundary.  \revone{There is no normal flow at the lower and upper boundaries.}
A range of mesh spacings are chosen so that $\Delta x \mathbin{:} \Delta z = 2\mathbin{:}1$ to match the original test specification from Sch\"{a}r et al. \citep{schaer2002}.

The terrain is wave-shaped, specified by the surface elevation $h$ such that
\begin{subequations}
\begin{align}
   h(x) &= h^\star \cos^2 ( \alpha x )
\shortintertext{where}
   h^\star(x) &= \left\{ \begin{array}{l l}
       h_0 \cos^2 ( \beta x ) & \quad \text{if $| x | < a$} \\
	0 & \quad \text{otherwise}
    \end{array} \right.
\end{align}
\end{subequations}
where $a = \SI{25}{\kilo\meter}$ is the mountain envelope half-width, $h_0 = \SI{6}{\kilo\meter}$ is the maximum mountain height, $\lambda = \SI{8}{\kilo\meter}$ is the wavelength, \(\alpha = \pi / \lambda\) and \(\beta = \pi / (2a)\).  Note that, in order to make this test more challenging, the mountain height $h_0$ is double the mountain height used by \citep{schaer2002}.

A basic terrain-following (BTF) mesh is constructed by using the terrain profile to modify the uniform mesh.
The BTF method uses a linear decay function so that mesh surfaces become horizontal at the top of the model domain \citep{galchen-somerville1975},
\begin{equation}
	z(x) = \left( H - h(x) \right) \left( z^\star / H \right) + h(x) \label{eqn:btf}
\end{equation}
where $z$ is the geometric height, $H$ is the height of the domain, $h(x)$ is the surface elevation and $z^\star$ is the computational height of a mesh surface.  If there were no terrain then $h = 0$ and $z = z^\star$.

A velocity field is prescribed with uniform horizontal flow aloft and zero flow near the ground,
\begin{align}
	u(z) = u_0 \left\{ \begin{array}{l l}
		1 & \quad \text{if $z \geq z_2$} \\
		\sin^2 \left( \frac{\pi}{2} \frac{z - z_1}{z_2 - z_1} \right) & \quad \text{if $z_1 < z < z_2$} \\
		0 & \quad \text{otherwise}
	\end{array} \right.	
\end{align}
where $u_0 = \SI{10}{\meter\per\second}$, $z_1 = \SI{7}{\kilo\meter}$ and $z_2 = \SI{8}{\kilo\meter}$.

A tracer with density $\phi$ has the shape
\begin{subequations}
\begin{align}
	\phi(x, z) &= \phi_0 \left\{ \begin{array}{l l}
		\cos^2 \left( \frac{\pi r}{2} \right) & \quad \text{if $r \leq 1$} \\
		0 & \quad \text{otherwise}
	\end{array} \right.
\intertext{with radius $r$ given by}
	r &= \sqrt{
		\left( \frac{x - x_0}{A_x} \right)^2 + 
		\left( \frac{z - z_0}{A_z} \right)^2
	}
\end{align}
\label{eqn:tracer}
\end{subequations}
where $A_x = \SI{25}{\kilo\meter}$, $A_z = \SI{3}{\kilo\meter}$ are the horizontal and vertical half-widths respectively, and $\phi_0 = \SI{1}{\kilogram\per\meter\cubed}$ is the maximum density of the tracer.  At $t = \SI{0}{\second}$, the tracer is centred at $(x_0, z_0) = (\SI{-50}{\kilo\meter}, \SI{12}{\kilo\meter})$ so that the tracer is upwind of the mountain, in the region of uniform flow above $z_2$.

\begin{figure}
	\centering
	\includegraphics{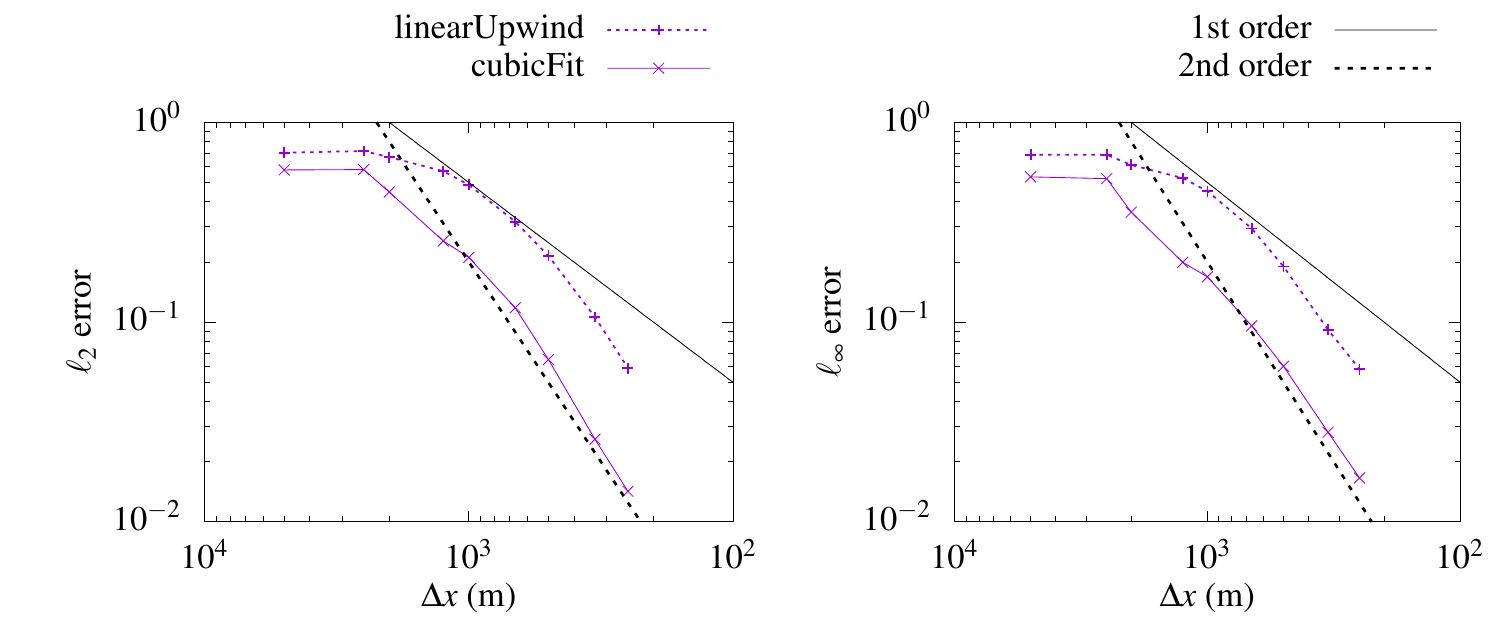}
	\caption{Numerical convergence of the two-dimensional horizontal transport test over mountains.  $\ell_2$ errors (equation~\ref{eqn:l2-error}) and $\ell_\infty$ errors (equation~\ref{eqn:linf-error}) are marked at mesh spacings between \SI{5000}{\meter} and \SI{250}{\meter} using linearUpwind and cubicFit transport schemes on basic terrain-following meshes.}
	\label{fig:horizontalAdvection-convergence}
\end{figure}

Tests are integrated for \SI{10000}{\second} using time-steps chosen for each mesh so that the maximum Courant number is about \num{0.4}.  \revone{This choice yields a time-step that is well below any stability limit, as recommended by Lauritzen et al. \citep{lauritzen2012}.}  By the end of integration the tracer is positioned downwind of the mountain.
The analytic solution at $t = \SI{10000}{\second}$ is centred at $(x_0, z_0) = (\SI{50}{\kilo\meter}, \SI{12}{\kilo\meter})$.  Error norms are calculated by subtracting the analytic solution from the numerical solution,
\begin{align}
	\ell_2 &= \sqrt{\frac{\sum_c \left(\phi - \phi_T \right)^2 \mathcal{V}_c}{\sum_c \left(\phi_T^2 \mathcal{V}_c \right)}} \label{eqn:l2-error} \\
	\ell_\infty &= \frac{\max_c |\phi - \phi_T|}{\max_c |\phi_T|} \label{eqn:linf-error}
\end{align}
where $\phi$ is the numerical value, $\phi_T$ is the analytic value, $\sum_c$ denotes a summation over all cells $c$ in the domain, and $\max_c$ denotes a maximum value of any cell.

Tests were performed using the linearUpwind and cubicFit schemes at mesh spacings between $\Delta x = \SI{250}{\meter}$ and $\Delta x = \SI{5000}{\meter}$.
Numerical convergence in the $\ell_2$ and $\ell_\infty$ norms is plotted in figure~\ref{fig:horizontalAdvection-convergence}.
The linearUpwind and cubicFit schemes are second-order convergent at all but the coarsest mesh spacings where errors are saturated for both schemes.

\revone{The cubicFit scheme achieves a given $\ell_2$ error using a mesh spacing that is almost twice as coarse as that needed by the linearUpwind scheme.  Doubling the mesh spacing results in a coarser mesh with four times fewer cells because the $\Delta x \mathbin{:} \Delta z$ aspect ratio is fixed.
Recall that the stencil for the cubicFit scheme has about twice as many cells as the stencil for the linearUpwind scheme.
Hence, for a given $\ell_2$ error, the computational cost during integration of the cubicFit scheme is about half the computational cost of the linearUpwind scheme.}

This test demonstrates that cubicFit is second-order convergent in the domain interior irrespective of mesh distortions.  The cubicFit scheme achieves In the next test, we assess the numerical accuracy of the transport schemes near a distorted, mountainous lower boundary.

%% file: mountainAdvection.tex
\subsection{Transport over a mountainous lower boundary}
\label{sec:mountainAdvection}

The horizontal transport test in the previous section is useful for assessing numerical accuracy on terrain-following meshes, but it presents no particular challenge on cut cell meshes because there is no flow through the distorted cut cells near the ground \citep{good2014}.
Here we present another variant of the standard horizontal transport test that challenges transport schemes on all mesh types.  By positioning the tracer next to the ground and modifying the velocity field, we can assess the accuracy of the cubicFit scheme near the lower boundary.  Results using the cubicFit scheme are compared with the linearUpwind scheme on basic terrain-following, cut cell and slanted cell meshes.

Cut cell meshes are constructed using the ASAM grid generator \citep{jaehn2015,asam2010}.  Slanted cell meshes are constructed following the approach by Shaw and Weller \citep{shaw-weller2016}: vertices that are underground are moved up to the surface and zero-area faces and zero-volume cells are removed.  Unlike \citep{shaw-weller2016}, vertices are never moved downwards.

\begin{figure}
	\centering
	\includegraphics{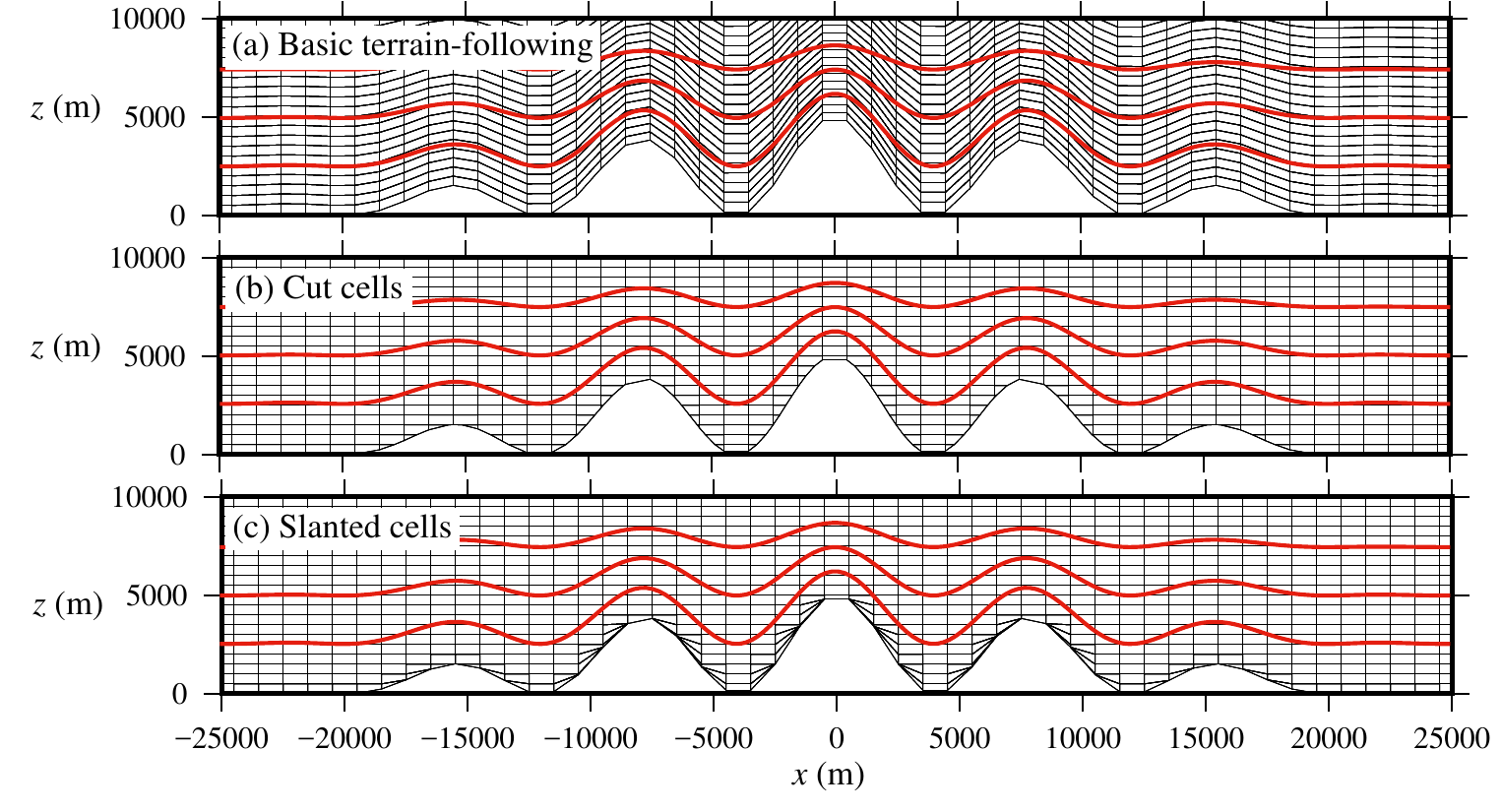}
	\caption{Two dimensional $x$-$z$ meshes created with the (a) basic terrain-following, (b) cut cell, and (c) slanted cell methods, and used for the tracer transport tests in section~\ref{sec:mountainAdvection}.  Cell edges are marked by thin black lines.  The peak mountain height $h_0 = \SI{5}{\kilo\meter}$.
The velocity field is the same for all mesh types with streamlines marked on each panel by thick red lines.  The velocity field (equation~\ref{eqn:streamfunc-btf}) follows the lower boundary and becomes entirely horizontal above $H_1 = \SI{10}{\kilo\meter}$.
Only the lowest \SI{10}{\kilo\meter} for the central region of the \revone{domain is shown}.  The entire domain is \SI{301}{\kilo\meter} wide and \SI{25}{\kilo\meter} high.}
	\label{fig:mountainAdvection-meshes}
\end{figure}

\revtwo{Following Sch\"ar et al. \citep{schaer2002}, the domain is \SI{301}{\kilo\meter} wide and \SI{25}{\kilo\meter} high as measured between parallel boundary edges, with a mesh spacing of $\Delta x = \SI{1000}{\meter}$ and $\Delta z = \SI{500}{\meter}$.}  \revother{The boundary conditions are the same as those used in section~\ref{sec:horizontal}.}
Cell edges in the central region of the domain are shown in figure~\ref{fig:mountainAdvection-meshes} for each of the three mesh types.
Cells in the BTF mesh are highly distorted over steep slopes (figure~\ref{fig:mountainAdvection-meshes}a) while the cut cell mesh (figure~\ref{fig:mountainAdvection-meshes}b) and slanted cell mesh (figure~\ref{fig:mountainAdvection-meshes}c) are orthogonal everywhere except for cells nearest the ground.

Similar to the approach by \citep{shaw-weller2016}, a velocity field is chosen that follows the terrain at the surface and becomes entirely horizontal aloft.
A streamfunction $\Psi$ is used so that the discrete velocity field is non-divergent, such that
\begin{equation}
	\Psi(x,z) = -u_0 H_1 \frac{z - h}{H_1 - h} \label{eqn:streamfunc-btf}
\end{equation}
where $u_0 = \SI{10}{\meter\per\second}$, which is the horizontal velocity where $h(x) = 0$.
There is no normal flow at \revone{the lower and upper boundaries}.  The velocity field becomes purely horizontal above $H_1 = \SI{10}{\kilo\meter}$ and, elsewhere, the flow is predominantely horizontal, with non-zero vertical velocities only above sloping terrain.
The horizontal and vertical components of velocity, $u$ and $w$, are given by
\begin{align}
	u &= -\frac{\partial \Psi}{\partial z} = u_0 \frac{H_1}{H_1 - h}, \quad w = \frac{\partial \Psi}{\partial x} = u_0 H_1 \frac{\mathrm{d} h}{\mathrm{d} x} \frac{H_1 - z}{\left( H_1 - h \right)^2} \label{eqn:uw-btf} \text{ ,}\\
	\frac{\mathrm{d} h}{\mathrm{d} x} &= - h_0 \left[ 
		\beta \cos^2 \left( \alpha x \right) \sin \left( 2 \beta x \right) +
		\alpha \cos^2 \left( \beta x \right) \sin \left( 2 \alpha x \right)
	\right] \text{ .}
\end{align}
Unlike the horizontal transport test in \citep{schaer2002}, the velocity field presented here extends from the top of the domain all the way to the ground.

The flow is deliberately misaligned with the BTF, cut cell and slanted cell meshes away from the ground (figure~\ref{fig:mountainAdvection-meshes}) to ensure that flow always crosses mesh surfaces in order to challenge the transport scheme.
The value of $H_1$ is chosen to be much smaller than the domain height $H$ in equation~\eqref{eqn:btf} so that flow crosses the surfaces of the BTF mesh.
This is evident in figure~\ref{fig:mountainAdvection-meshes}a where the the velocity streamlines are tangential to the mesh only at the ground.

\begin{table}
	\centering
\begin{tabular}{l S S S S S}
\hline
	& \multicolumn{5}{c}{Peak mountain height $h_0$ (\si{\kilo\meter})} \\
	Mesh type & 0 & 3 & 4 & 5 & 6 \\
\hline
	BTF & 40 & 16 & 10 & 8 & 5 \\
	Cut cell & 40 & 1.6 & 1.6 & 0.5 & 1.5  \\
	Slanted cell & 40 & 8 & 6.25 & 5 & 4  \\
\hline
\end{tabular}
	\caption{Time-steps (\si{\second}) for the two-dimensional transport test over a mountainous lower boundary.  The time-steps were chosen so that the maximum Courant number was between \num{0.36} and \num{0.46}.}
	\label{tab:mountainAdvection:timesteps}
\end{table}

The tracer is again defined by equation~\eqref{eqn:tracer} but is now positioned at the ground with $(x_0, z_0) = (\SI{-50}{\kilo\meter}, \SI{0}{\kilo\meter})$ with half-widths $A_x = \SI{25}{\kilo\meter}$ and $A_z = \SI{10}{\kilo\meter}$.
Tests are integrated forward for \SI{10000}{\second}.  The time-step was chosen for each mesh so that the maximum Courant number was about \num{0.4} (table~\ref{tab:mountainAdvection:timesteps}).
An analytic solution at \SI{10000}{\second} is obtained by calculating the new horizontal position of the tracer.  Integrating along the trajectory yields $t$, the time taken to move from the left side of the mountain to the right \citep{shaw-weller2016}:
\begin{align}
	\mathrm{d}t &= \mathrm{d}x / u(x) \\
	t &= \int_0^x \frac{H - h(x)}{u_0 H}\:\mathrm{d}x \\
	t &= \frac{x}{u_0} - \frac{h_0}{16 u_0 H} \left[ 4x + \frac{\sin 2 (\alpha + \beta) x}{\alpha + \beta} 
 \frac{\sin 2(\alpha - \beta) x}{\alpha - \beta} + 2 \left( \frac{\sin 2\alpha x}{\alpha} + \frac{\sin 2\beta x}{\beta} \right) \right]
\end{align}
By solving this equation we find that \(x(t=\SI{10000}{\second}) = \SI{54342.8}{\meter}\) when $h_0 = \SI{5}{\kilo\meter}$.

The tracer density boundary conditions are the same as those in the previous test.
Since the cubicFit transport scheme uses values at boundaries with Dirichlet boundary conditions, the cubicFit scheme uses only inlet boundary values in this test case.

\begin{figure}
	\centering
	\includegraphics{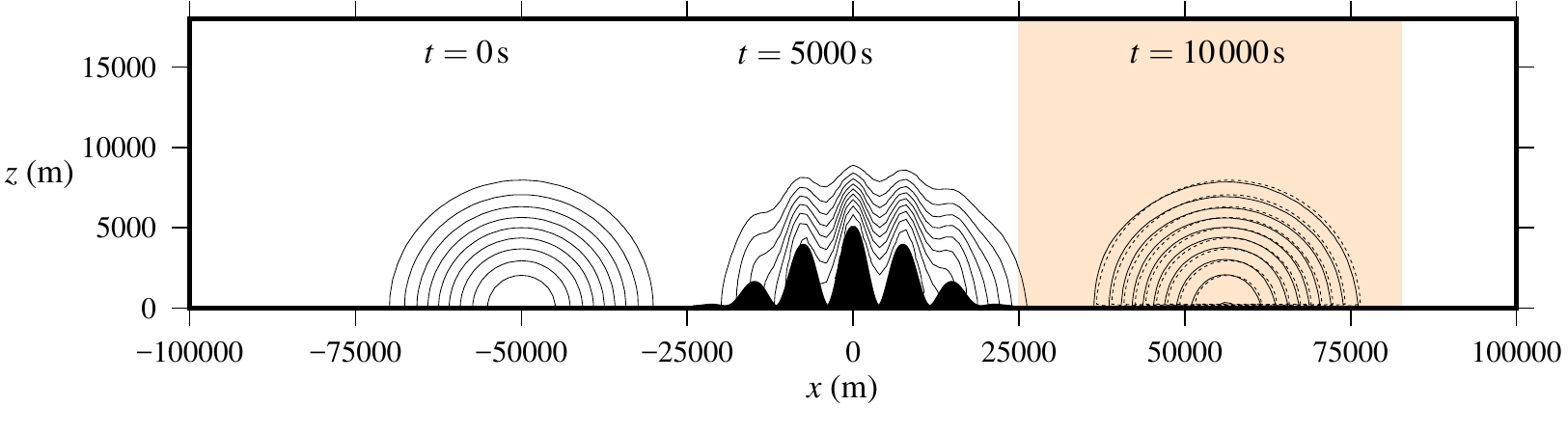}
	\caption{Evolution of the tracer in the two-dimensional transport test over a mountainous lower boundary.  The tracer is transported to the right over the wave-shaped terrain.  Tracer contours are every \SI{0.1}{\kilo\gram\per\meter\cubed}.  The result obtained using the cubicFit scheme on the basic terrain-following mesh is shown at $t=\SI{0}{\second}$, $t=\SI{5000}{\second}$ and $t=\SI{10000}{\second}$ with solid black contours. The analytic solution at $t=\SI{10000}{\second}$ is shown with dotted contours.
	The shaded box indicates the region that is plotted in figure~\ref{fig:mountainAdvection-errors}.}
	\label{fig:mountainAdvection-tracer}
\end{figure}

Three series of tests were performed using similar configurations.  The first series uses a peak mountain height of $h_0 = \SI{5}{\kilo\meter}$ to examine errors on different mesh types using the two transport schemes.
The second series varies the peak mountain height to examine the sensitivity of the transport schemes to mesh distortions.
The third series verifies accuracy at Courant numbers close to \revother{the limit of stability}, and examines the longest stable time-step for different mesh types.

For the first series of tests with $h_0 = \SI{5}{\kilo\meter}$, tracer contours at the initial time $t=\SI{0}{\second}$, half-way time $t=\SI{5000}{\second}$, and end time $t=\SI{10000}{\second}$ are shown in figure~\ref{fig:mountainAdvection-tracer} using the cubicFit scheme on the BTF mesh.  As apparent at $t=\SI{5000}{\second}$, the tracer is distorted by the terrain-following velocity field as it passes over the mountain, but its original shape is restored once it has cleared the mountain by $t=\SI{10000}{\second}$.
A small phase lag is apparent when the numerical solution marked with solid contour lines is compared with the analytic solution marked with dotted contour lines.

\begin{figure}
	\centering
	\includegraphics{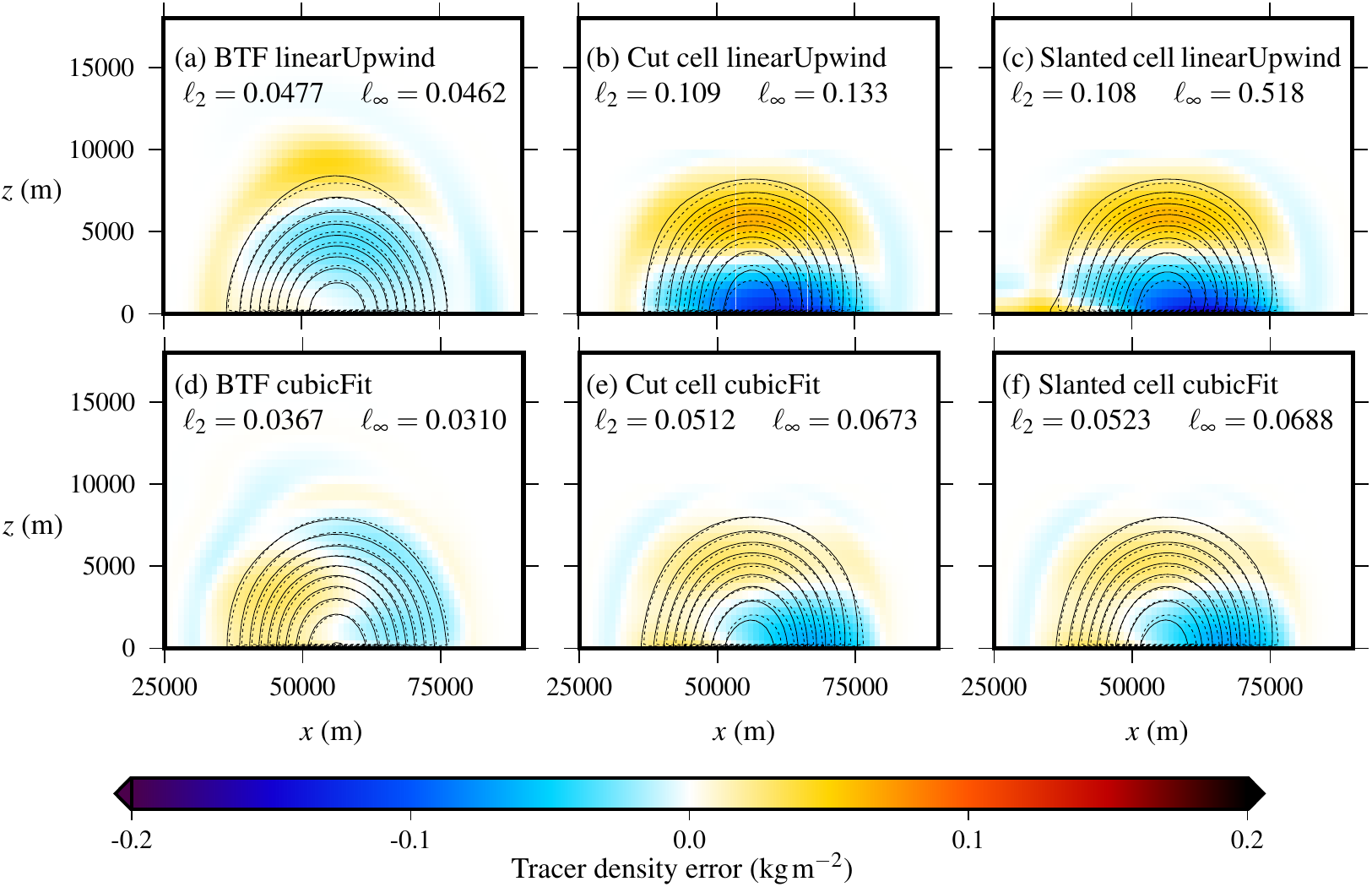}
	\caption{Tracer contours at $t=\SI{10000}{\second}$ for the two-dimensional tracer transport tests over a mountainous lower boundary.  A region in the lee of the mountain is plotted corresponding to the shaded area in figure~\ref{fig:mountainAdvection-tracer}.  Results are presented on BTF, cut cell and slanted cell meshes (shown in figure~\ref{fig:mountainAdvection-meshes}) using the linearUpwind and cubicFit transport schemes.  The numerical solutions are marked by solid black lines.  The analytic solution is marked by dotted lines.  Contours are every \SI{0.1}{\kilo\gram\per\meter\cubed}.}
	\label{fig:mountainAdvection-errors}
\end{figure}

Numerical errors are more clearly revealed by subtracting the analytic solution from the numerical solution.
Error fields are compared between BTF, cut cell and slanted cell meshes using the linearUpwind scheme (figures~\ref{fig:mountainAdvection-errors}a, \ref{fig:mountainAdvection-errors}b and \ref{fig:mountainAdvection-errors}c respectively) and the cubicFit scheme (figures~\ref{fig:mountainAdvection-errors}d, \ref{fig:mountainAdvection-errors}e and \ref{fig:mountainAdvection-errors}f respectively).
Results are least accurate using the linearUpwind scheme on the slanted cell mesh (figure~\ref{fig:mountainAdvection-errors}c).  The final tracer is slightly distorted and does not extend far enough towards the ground.
The $\ell_\infty$ error magnitude is reduced by using the linearUpwind scheme on the cut cell mesh (figure~\ref{fig:mountainAdvection-errors}b), but the shape of the error remains the same.
The cubicFit scheme is less sensitive to the choice of mesh with similar error magnitudes on the BTF mesh (figure~\ref{fig:mountainAdvection-errors}d), cut cell mesh (figure~\ref{fig:mountainAdvection-errors}e) and slanted cell mesh (figure~\ref{fig:mountainAdvection-errors}f).  Errors using the cubicFit scheme on cut cell and slanted cell meshes are much smaller than the errors using the linearUpwind scheme on the same meshes.

\begin{figure}
	\centering
	\includegraphics{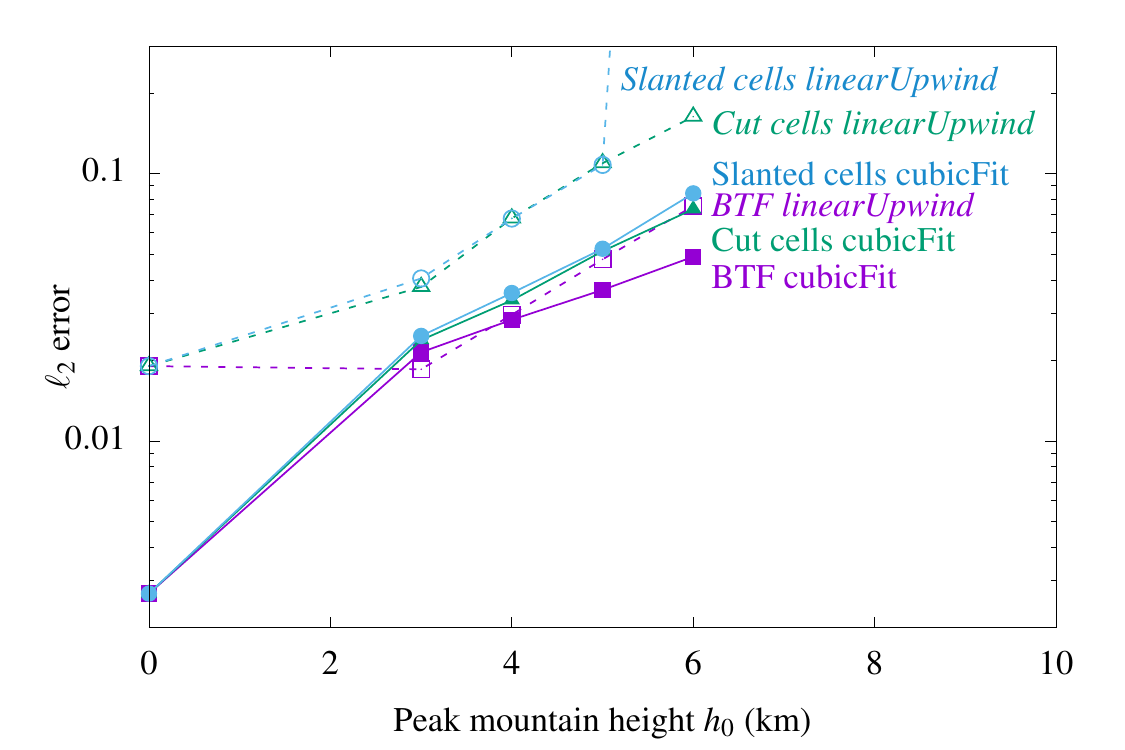}
	\caption{Error measures for the two-dimensional tracer transport tests over a mountainous lower boundary.  Peak mountain heights $h_0$ are from \SIrange{0}{6}{\kilo\meter}.  Results are compared on BTF, cut cell and slanted cell meshes using the linearUpwind and the cubicFit schemes.  At $h_0 = \SI{0}{\kilo\meter}$ the terrain is entirely flat and the BTF, cut cell and slanted cell meshes are identical.  At $h_0 = \SI{6}{\kilo\meter}$ the linearUpwind scheme is unstable on the slanted cell mesh.}
	\label{fig:mountainAdvection-l2ByMountainHeight}
\end{figure}

To further examine the performance of the cubicFit scheme in the presence of steep terrain, a second series of tests were performed in which the peak mountain height was varied from \SIrange{0}{6}{\kilo\meter} keeping all other parameters constant.
Results were obtained on BTF, cut cell and slanted cell meshes using the linearUpwind scheme and cubicFit scheme.  Again, the time-step was chosen for each test so that the maximum Courant number was about \num{0.4} (table~\ref{tab:mountainAdvection:timesteps}).  The $\ell_2$ error was calculated by subtracting the analytic solution from the numerical solution (figure~\ref{fig:mountainAdvection-l2ByMountainHeight}).
Note that the analytic solution is a function of mountain height, with the tracer travelling farther over higher mountains due to non-divergent flow through a narrower channel.
In all cases, error increases with increasing mountain height because steeper slopes lead to greater mesh distortions.
Errors are identical for a given transport scheme when $h_0 = \SI{0}{\kilo\meter}$ and the ground is entirely flat because the BTF, cut cell and slanted cell meshes are identical.
Compared with the cubicFit scheme, the linearUpwind scheme is more sensitive to the mesh type and mountain height.  The linearUpwind scheme is unstable on the slanted cell mesh with a peak mountain height $h_0 = \SI{6}{\kilo\meter}$ despite using a Courant number of \num{0.428}.
In contrast, the cubicFit scheme is less sensitive to the mesh type and errors grow more slowly with increasing mountain height.  The cubicFit scheme yields stable results in all tests.

\begin{figure}
	\centering
	\includegraphics{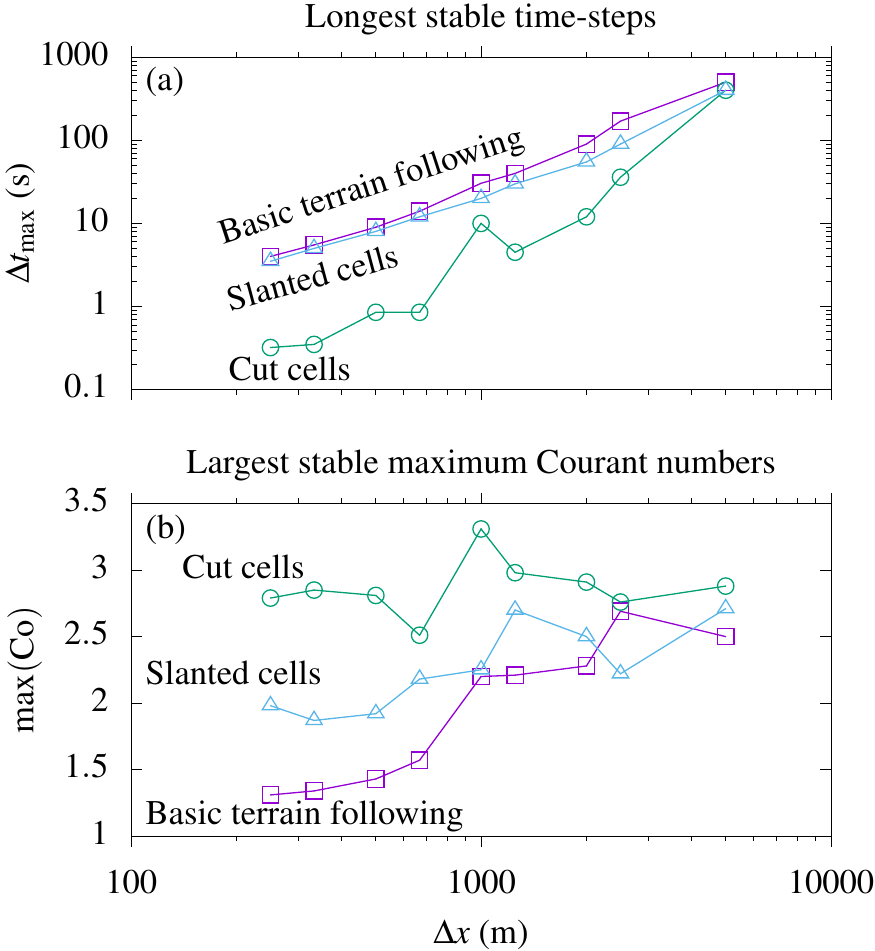}
	\caption{\revtwo{(a) Longest stable time-steps, $\Delta t_\mathrm{max}$, and (b) largest stable maximum Courant numbers, $\max(\mathrm{Co})$, for the two-dimensional tracer transport test over a mountainous lower boundary.  Results were obtained on basic terrain-following, cut cell and slanted cell meshes at mesh spacings between $\Delta x = \SI{5000}{\meter}$ and $\Delta x = \SI{250}{\meter}$.  The largest stable maximum Courant numbers were calculated from the corresponding longest stable time-steps using equation~\eqref{eqn:co}.}}
	\label{fig:mountainAdvection-maxdt}
\end{figure}

A final series of tests were performed \revtwo{to determine the stability limit of the cubicFit scheme with the two-stage Heun time-stepping scheme (equation~\ref{eqn:heun}).}
The tracer was transported on BTF, slanted cell and cut cell meshes with a variety of mesh spacings between $\Delta x = \SI{5000}{\meter}$ and $\Delta x = \SI{125}{\meter}$.  $\Delta z$ was chosen so that a constant aspect ratio is preserved such that $\Delta x / \Delta z = 2$.
\revtwo{For each test, the time-step was increased until the result became unstable.  The largest stable time-steps, $\Delta t_\mathrm{max}$, are presented in figure~\ref{fig:mountainAdvection-maxdt}a.}
BTF meshes permit the longest time-steps of all three mesh types since cells are almost uniform in volume.  As expected, the longest stable time-step scales linearly with BTF mesh spacing.
There is no such linear scaling on cut cell meshes because these meshes can have arbitrarily small cells.  The time-step constraints on cut cell meshes are the most severe of the three mesh types.  Slanted cell meshes have a slightly stronger time-step constraint than BTF meshes but still exhibit similar linear scaling with mesh spacing.  Furthermore, a dynamical model that uses slanted cell meshes instead of BTF meshes is expected to calculate pressure gradients more accurately \citep{shaw-weller2016}.

\revtwo{Figure~\ref{fig:mountainAdvection-maxdt}b presents the largest stable maximum Courant numbers, $\max(\mathrm{Co})$, which were calculated by substituting $\Delta t = \Delta t_\mathrm{max}$ into equation~\eqref{eqn:co}.
On basic terrain following meshes, the maximum Courant number tends towards about \num{1.3} with finer mesh spacings.
No such trend is found on cut cell or slanted cell meshes.
Cut cell meshes permit the largest maximum Courant numbers of around \num{3}, but the largest stable time-steps on cut cell meshes are still smaller than corresponding time-steps on basic terrain following and slanted cell meshes.}

\revtwo{This paper focuses on the spatial discretisation of the cubicFit scheme, but the stability limit depends also upon the choice of time-stepping.  As such, we have not calculated a theoretical Courant number limit, although such an analysis should be possible using the techniques in \citep{baldauf2008}.}

The transport tests presented in this section demonstrate that the cubicFit scheme is suitable for flows over very steep terrain on two-dimensional terrain-following, cut cell and slanted cell meshes.  The cubicFit scheme is less sensitive to the mesh type and mountain steepness compared to the linearUpwind scheme.  The linearUpwind scheme becomes unstable over very steep slopes but the cubicFit scheme is stable for all tests.
The \revother{accuracy of the} \revtwo{cubicFit scheme} was largely insensitive to the choice of time-step.
In the next section, we evaluate the cubicFit scheme using more complex, deformational flows on icosahedral meshes and cubed-sphere meshes.

%% file: deformationSphere.tex
\subsection{Deformational flow on a sphere}
\label{sec:deformationSphere}
The tests so far have used flows that are mostly uniform on meshes that are based on rectangular cells.
To ensure that the cubicFit transport scheme is suitable for complex flows on a variety of meshes, we use a standard test of deformational flow on a spherical Earth, as specified by Lauritzen et al. \citep{lauritzen2012}.  
Results are compared between linearUpwind and cubicFit schemes using \revone{hexagonal-icosahedral meshes} and cubed-sphere \revother{meshes}.
Hexagonal-icosahedral meshes are constructed by successive refinement of a regular icosahedron following the approach by \revtwo{\citep{thuburn2014,heikes-randall1995a,heikes-randall1995b} without any mesh twisting}.
Cubed-sphere meshes are constructed using an equi-distant gnomic projection of a cube having a uniform Cartesian mesh on each panel \citep{staniforth-thuburn2012}.

Following appendix A9 in \citep{lauritzen2014}, the average equatorial spacing $\Delta \lambda$ is used as a measure of mesh spacing.  It is defined as
\begin{align}
	\Delta \lambda = \ang{360} \frac{\overline{\Delta x}}{2 \pi R_e}
\end{align}
where $\overline{\Delta x}$ is the mean distance between cell centres and $R_e = \SI{6.3712e6}{\meter}$ is the radius of the Earth.

The deformational flow test specified by Lauritzen et al. \citep{lauritzen2012} comprised six elements:
\begin{enumerate}
\item a convergence test using a Gaussian-shaped tracer
\item a ``minimal'' resolution test using a cosine bell-shaped tracer
\item a test of filament preservation
\item a test using a ``rough'' slotted cylinder tracer
\item a test of correlation preservation between two tracers
\item a test using a divergent velocity field
\end{enumerate}
We assess the cubicFit scheme using the first two tests only.  We do not consider filament preservation, correlation preservation, or the transport of a ``rough'' slotted cylinder because no shape-preserving filter has yet been developed for the cubicFit scheme.  Stable results were obtained when testing the cubicFit scheme using a divergent velocity field, but no further analysis is made here.

The first deformational flow test uses an infinitely continuous initial tracer that is transported in a non-divergent, time-varying, rotational velocity field.
The velocity field deforms two Gaussian `hills' of tracer into thin vortical filaments.  Half-way through the integration the rotation reverses so that the filaments become circular hills once again.  The analytic solution at the end of integration is identical to the initial condition.
A rotational flow is superimposed on a time-invariant background flow in order to avoid error cancellation.
The non-divergent velocity field is defined by the streamfunction $\Psi$,
\begin{align}
	\Psi(\lambda, \theta, t) = \frac{10 R_e}{T} \sin^2 \left(\lambda'\right) \cos^2 \left(\theta\right) \cos \left( \frac{\pi t}{T} \right) - \frac{2 \pi R_e}{T} \sin\left(\theta\right)
\end{align}
where $\lambda$ is a longitude, $\theta$ is a latitude, $\lambda' = \lambda - 2 \pi t / T$, and $T = \num{12}$ days is the duration of integration.  The time-step is chosen such that the maximum Courant number is about 0.4.

The initial tracer density $\phi$ is defined as the sum of two Gaussian hills,
\begin{align}
	\phi = \phi_1(\lambda, \theta) + \phi_2(\lambda, \theta) \text{ .}
\end{align}
An individual hill $\phi_i$ is given by
\begin{align}
	\phi_i(\lambda, \theta) = \phi_0 \exp\left( -b \left( \frac{|\mathbf{x} - \mathbf{x}_i|}{R_e} \right)^2 \right)
\end{align}
where $\phi_0 = \SI{0.95}{\kilo\gram\per\meter\cubed}$ and $b = 5$.  The Cartesian position vector $\mathbf{x} = (x,y,z)$ is related to the spherical coordinates $(\lambda, \theta)$ by
\begin{align}
	(x,y,z) = (R_e \cos \theta \cos \lambda, R_e \cos \theta \sin \lambda, R_e \sin \theta) \label{eqn:spherical-cartesian} \text{ .}
\end{align}
The centre of hill $i$ is positioned at $\mathbf{x}_i$.  In spherical coordinates, two hills are centred at
\begin{align}
	(\lambda_1,\theta_1) &= (5 \pi /6, 0) \\
	(\lambda_2,\theta_2) &= (7 \pi /6, 0)
\end{align}

\begin{figure}
	\centering
	\includegraphics{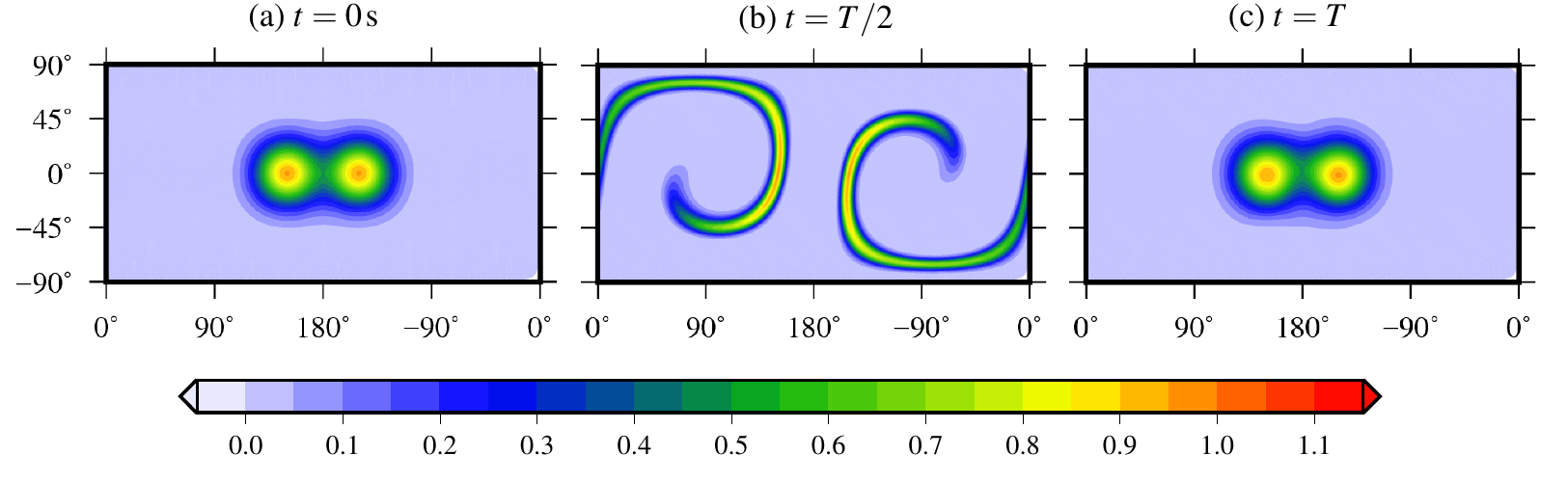}
	\caption{Tracer fields for the deformational flow test using initial Gaussian hills.  The tracer is deformed by the velocity field before the rotation reverses to return the tracer to its original distribution: (a) the initial tracer distribution at $t = \SI{0}{\second}$; (b) by $t=T/2$ the Gaussian hills are stretched into a thin S-shaped filament; (c) at $t=T$ the tracer resembles the initial Gaussian hills except for some distortion and diffusion due to numerical errors.  Results were obtained with the cubicFit scheme on a \revone{hexagonal-icosahedral mesh} with an average equatorial mesh spacing of $\Delta \lambda = \ang{0.542}$.}
	\label{fig:deformationSphere-evolution}
\end{figure}

The results in figure~\ref{fig:deformationSphere-evolution} are obtained using the cubicFit scheme on a hexagonal-icosahedral mesh with $\Delta \lambda = \ang{0.542}$.  The initial Gaussian hills are shown in figure~\ref{fig:deformationSphere-evolution}a.  At $t=T/2$ the tracer has been deformed into an S-shaped filament (figure~\ref{fig:deformationSphere-evolution}b).  By $t=T$ the tracer has almost returned to its original distribution except for some slight distortion and diffusion that are the result of numerical errors (figure~\ref{fig:deformationSphere-evolution}c).

\begin{figure}
	\centering
	\includegraphics{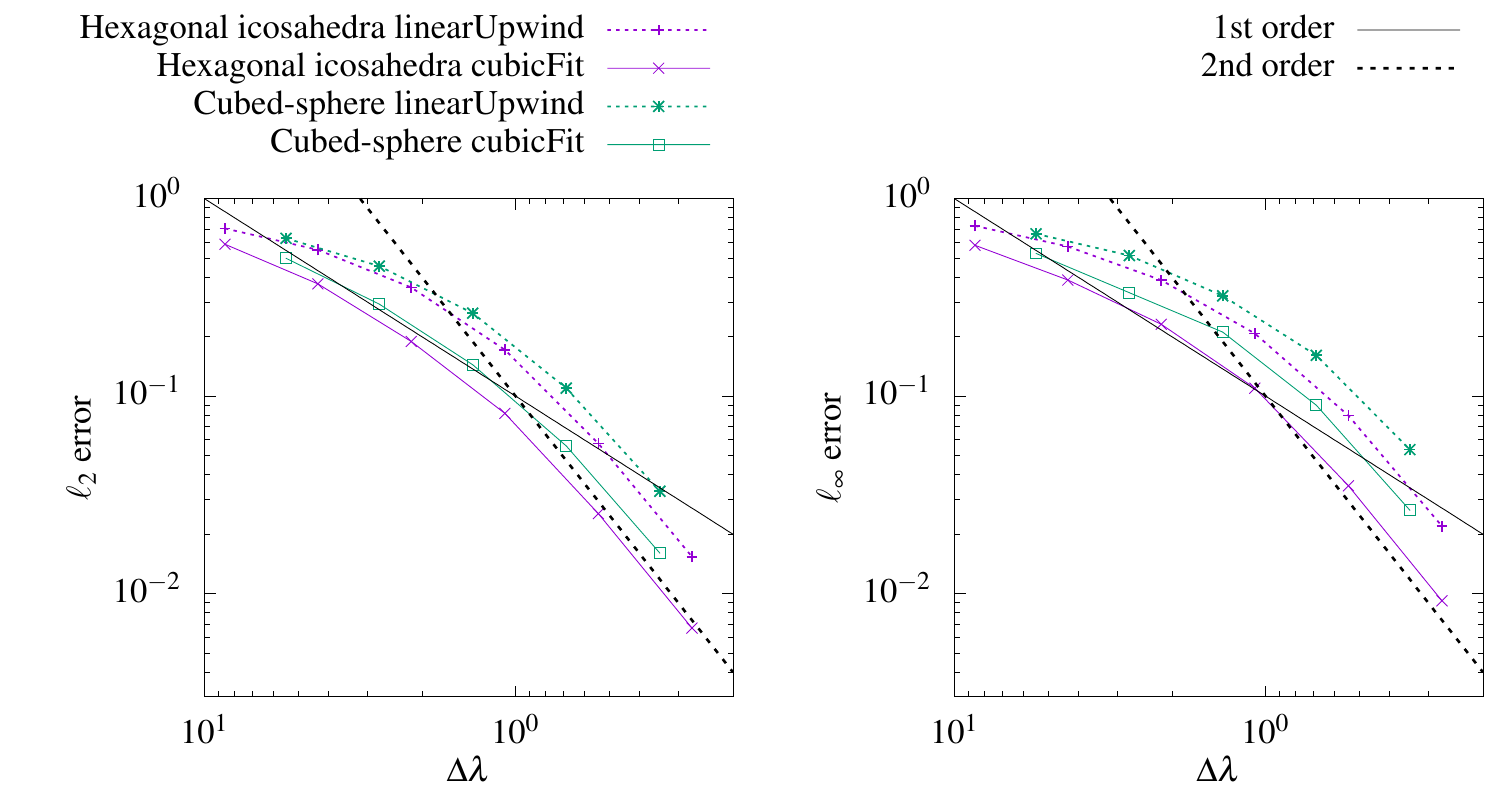}
	\caption{Numerical convergence of the deformational flow test on the sphere using initial Gaussian hills.  $\ell_2$ errors (equation~\ref{eqn:l2-error}) and $\ell_\infty$ errors (equation~\ref{eqn:linf-error}) are marked at mesh spacings between \ang{8.61} and \ang{0.271} using the linearUpwind scheme (dotted lines) and the cubicFit scheme (solid lines) on \revone{hexagonal-icosahedral meshes} and cubed-sphere meshes.}
	\label{fig:deformationSphere-gaussian-convergence}
\end{figure}

To determine the order of convergence and relative accuracy of the linearUpwind and cubicFit schemes, the same test was performed at a variety of mesh spacings betweeen $\Delta \lambda = \ang{8.61}$ and $\Delta \lambda = \ang{0.271}$ on \revone{hexagonal-icosahedral meshes} and cubed-sphere meshes.  The results are shown in figure~\ref{fig:deformationSphere-gaussian-convergence}.
The solution is slow to converge at coarse resolutions, and this behaviour agrees with the results from Lauritzen et al. \citep{lauritzen2012}.  Both linearUpwind and cubicFit schemes achieve second-order accuracy at smaller mesh spacings. 
For any given mesh type and mesh spacing, the cubicFit scheme is more accurate than the linearUpwind scheme.
Results are more accurate using \revone{hexagonal-icosahedral meshes} compared to cubed-sphere meshes.  It is not known whether the larger errors on cubed-sphere meshes are due to mesh non-uniformities at panel corners but there is no evidence of grid \revone{imprinting in} the error fields (not shown).

A slightly more challenging variant of the same test is performed using a quasi-smooth tracer field defined as the sum of two cosine bells,
\begin{align}
	\phi =
	\begin{cases}
		b + c \phi_1(\lambda, \theta) & \quad \text{if $r_1 < r$,} \\
		b + c \phi_2(\lambda, \theta) & \quad \text{if $r_2 < r$,} \\
		b			      & \quad \text{otherwise.}
	\end{cases}
\end{align}
The velocity field is the same as before.  This test is used to determine the ``minimal'' resolution, $\Delta \lambda_m$, which is specified by Lauritzen et al. \citep{lauritzen2012} as the coarsest mesh spacing for which $\ell_2 \approx 0.033$.

\begin{table}
	\robustify\it
	\centering
	\begin{tabular}{l l S[detect-weight, detect-shape, detect-mode]}
\hline
	Transport scheme & Mesh type & {Minimal resolution (\si{\degree})} \\
\hline
	linearUpwind & Cubed-sphere & \it 0.15 \\
	FARSIGHT, grid-point semi-Lagrangian \citep{white-dongarra2011} & Cubed-sphere & 0.1875 \\
	linearUpwind & \revone{Hexagonal-icosahedral} & \it 0.2 \\
	SLFV-SL, swept-area scheme \citep{miura2007} & \revone{Hexagonal-icosahedral} & 0.25 \\
	cubicFit & Cubed-sphere & \it 0.25 \\
	cubicFit & \revone{Hexagonal-icosahedral} & 0.3 \\
	ICON-FFSL, swept-area scheme \citep{miura2007} & \revone{Triangular-icosahedral} & 0.42 \\
\hline
\end{tabular}
	\caption{Minimal resolutions for the cubicFit and linearUpwind schemes in the test of deformational flow using cosine bells.  Italicised values have been extrapolated using the second-order convergence obtained at coarser mesh spacings.  For comparison with existing models, some results are also included for unlimited versions of the transport schemes from the intercomparison by Lauritzen et al. \citep{lauritzen2014}.}
	\label{tab:deformationSphere:minimal-resolution}
\end{table}

The minimal resolution for the cubicFit scheme on a hexagonal-icosahedral mesh is about $\Delta \lambda_m = \ang{0.3}$.  Tests were not performed at mesh spacings finer than $\Delta \lambda = \ang{0.271}$ but approximate minimal resolutions have been extrapolated from the second-order convergence that is found at fine mesh spacings.  These minimal resolutions are presented in table~\ref{tab:deformationSphere:minimal-resolution} along with a selection of transport schemes having similar minimal resolutions from the model intercomparison by Lauritzen et al. \citep{lauritzen2014}.

The series of deformational flow tests presented here demonstrate that the cubicFit scheme is suitable for transport on spherical meshes based on quadrilaterals and hexagons.  The cubicFit scheme is largely insensitive to the mesh type, and results are more accurate compared to the linearUpwind scheme for a given mesh type and mesh spacing.  Neither scheme requires special treatment at the corners of cubed-sphere panels.

%% file: conclusion.tex
\section{Conclusion}
\label{sec:conclusion}

Atmospheric models are using increasingly fine horizontal mesh spacings that resolve steep slopes in terrain resulting in highly-distorted meshes, increased numerical errors and numerical instabilities.
We have presented a new multidimensional method-of-lines transport scheme, cubicFit, that applies constraints derived from a von Neumann stability analysis to make the scheme stable over steep terrain on \revone{highly-distorted, arbitrary meshes}.
The scheme has a low computational cost at runtime, requiring only $n$ multiplies per face per \revtwo{time-stage} using a stencil with $n$ cells.  Stability constraint calculations are pre-computed during model initialisation since they depend upon the mesh geometry only.

The cubicFit scheme was compared to a multidimensional linear upwind scheme using three idealised numerical tests.
The first test transported a tracer horizontally above steep slopes on highly-distorted, two-dimensional terrain-following meshes.  The cubicFit scheme was second-order convergent regardless of mesh distortions.
The second test transported a tracer over a mountainous lower boundary using terrain-following, cut cell and slanted cell meshes.
The cubicFit scheme was generally insensitive to the type of mesh and less sensitive to terrain steepness compared to the multidimensional linear upwind scheme, and the scheme maintained accuracy up to its stability limit.
The third test evaluated the transport schemes in a standard deformational flow field on \revone{hexagonal-icosahedral meshes} and cubed-sphere meshes.
In all tests, compared to the multidimensional linear upwind scheme, the cubicFit transport scheme was more stable and more accurate.

%% file: vonNeumann.tex
\section*{Appendix A: One-dimensional von Neumann stability analysis}
Two analyses are performed in order to find stability constraints on the weights $\mathbf{w} = \mathbf{\tilde{b}_1^+} \cdot \mathbf{m}$ as appear in equation~\eqref{eqn:weightedPinv}.  The first analysis uses a two-cell approximation to derive separate constraints on the upwind weight $w_u$ and downwind weight $w_d$.  The second analysis uses three cells to derive a constraint that considers all weights in a stencil.

\subsection*{Two-cell analysis}
We start with the conservation equation for a dependent variable $\phi$ that is discrete-in-space and continuous-in-time
\begin{align}
\frac{\partial \phi_j}{\partial t} &= - v \frac{\phi_R - \phi_L}{\Delta x} \label{eqn:advectionLR} \\
\intertext{where $v$ is the velocity, and the left and right fluxes, $\phi_L$ and $\phi_R$, are weighted averages of the neighbouring cell centres.  Assuming that $v$ is positive}
\phi_L &= \alpha_u \phi_{j-1} + \alpha_d \phi_j \\
\phi_R &= \beta_u \phi_j + \beta_d \phi_{j+1}
\end{align}
where $\phi_{j-1}, \phi_j, \phi_{j+1}$ are cell centre values, and $j$ denotes a cell centre position $x = j \Delta x$ where $\Delta x$ is a uniform mesh spacing.
$\alpha_u$ and $\beta_u$ are the upwind weights and $\alpha_d$ and $\beta_d$ are the downwind weights for the left and right fluxes respectively, and $\alpha_u + \alpha_d = 1$ and $\beta_u + \beta_d = 1$.

At a given time $t = n \Delta t$ at time-level $n$ and with a time-step $\Delta t$, we assume a wave-like solution with an amplification factor $A$, such that
\begin{align}
	\phi_j^{(n)} &= A^n e^{\iu j k \Delta x} \label{eqn:vn}
\end{align}
where $\phi_j^{(n)}$ denotes a value of $\phi$ at position $j$ and time-level $n$.  Using this to rewrite the left-hand side of equation~\eqref{eqn:advectionLR}
\begin{align}
\frac{\partial \phi_j}{\partial t} &= \frac{\partial}{\partial t} \left( A^{t / \Delta t} \right) e^{ijk\Delta x} = \frac{\ln A}{\Delta t} A^n e^{ikj\Delta x} \\
\shortintertext{hence equation~\eqref{eqn:advectionLR} becomes}
\frac{\ln A}{\Delta t} &= - \frac{v}{\Delta x} \left( \beta_u + \beta_d e^{ik\Delta x} - \alpha_u e^{-ik\Delta x} - \alpha_d \right) \\
\ln A &= -c \left( \beta_u - \alpha_d + \beta_d \cos k\Delta x + \iu \beta_d \sin k \Delta x - \alpha_u \cos k\Delta x + \iu \alpha_u \sin k\Delta x \right)
\intertext{where the Courant number $c = v \Delta t / \Delta x$.
Let $\Re = \beta_u - \alpha_d + \beta_d \cos k\Delta x - \alpha_u \cos k\Delta x$ and
$\Im = \beta_d \sin k \Delta x + \alpha_u \sin k\Delta x$, then}
\ln A &= -c \left( \Re + \iu \Im \right) \\
A &= e^{-c \Re} e^{-\iu c \Im} \\
\shortintertext{and the complex modulus of $A$ is}
|A| &= e^{-c \Re} = \exp \left( -c \left( \beta_u - \alpha_d + \left(\beta_d - \alpha_u \right) \cos k\Delta x \right) \right) \text{ .}
\end{align}
For stability we need $|A| \leq 1$ and, imposing the additional constraints that $\alpha_u = \beta_u$ and $\alpha_d = \beta_d$, then
\begin{align}
\left( \alpha_u - \alpha_d \right) \left( 1 - \cos k\Delta x \right) &\geq 0 \quad \forall k\Delta x
\shortintertext{and, given $0 \leq 1 - \cos k \Delta x \leq 2$, then}
\alpha_u - \alpha_d \geq 0 \text{ .} \label{eqn:twopoint-lower}
\end{align}
Additionally, we do not want more damping than a first-order upwind scheme (where $\alpha_u = \beta_u = 1$, $\alpha_d = \beta_d = 0$), having an amplification factor, $A_\mathrm{up}$, so we need $\left\lvert A \right\rvert \geq \left\lvert A_\mathrm{up} \right\rvert$, hence
\begin{align}
	\exp \left( -c \left(\alpha_u - \alpha_d\right) \left( 1 - \cos k\Delta x \right) \right) &\geq \exp \left( -c \left(1 - \cos k\Delta x \right) \right) \quad \forall k\Delta x
\shortintertext{therefore}
	\alpha_u - \alpha_d &\leq 1 \text{ .} \label{eqn:twopoint-upper}
\end{align}
Now, knowing that $\alpha_u + \alpha_d = 1$ (or $\alpha_d = 1 - \alpha_u$) then, using equations~\eqref{eqn:twopoint-lower} and \eqref{eqn:twopoint-upper},
\begin{align}
	0.5 \leq \alpha_u &\leq 1 \text{ and} \label{eqn:vn:upwind} \\
	0 \leq \alpha_d &\leq 0.5 \label{eqn:vn:downwind} \text{ .}
\end{align}

\subsection*{Three-cell analysis}
We start again from equation~\eqref{eqn:advectionLR} but this time approximate $\phi_L$ and $\phi_R$ using three cell centre values,
\begin{align}
	\phi_L &= \alpha_{uu} \phi_{j-2} + \alpha_u \phi_{j-1} + \alpha_d \phi_j \\
	\phi_R &= \alpha_{uu} \phi_{j-1} + \alpha_u \phi_j + \alpha_d \phi_{j+1}
\end{align}
having used the same weights $\alpha_{uu}$, $\alpha_u$ and $\alpha_d$ for both left and right fluxes.
Substituting equation~\eqref{eqn:vn} into equation~\eqref{eqn:advectionLR} we find
\begin{align}
A = \exp\left( -c \left[ \alpha_{uu} \left( e^{-ik\Delta x} - e^{-2ik\Delta x} \right) + \alpha_u \left( 1 - e^{-ik\Delta x} \right) + \alpha_d \left( e^{ik\Delta x} - 1 \right) \right] \right)
\intertext{so that, if the complex modulus $\left\lvert A \right\rvert \leq 1$ then}
\alpha_u - \alpha_d + \left( \alpha_{uu} - \alpha_u + \alpha_d \right) \cos k\Delta x - \alpha_{uu} \cos 2k\Delta x \geq 0 \text{ .}
\end{align}
If $k\Delta x = \pi$ then $\cos k\Delta x = -1$ and $\cos 2k\Delta x = 1$ and $\alpha_u - \alpha_d \geq \alpha_{uu}$.  If $k\Delta x = \pi / 2$ then $\cos k\Delta x = 0$ and $\cos 2k\Delta x = -1$ and $\alpha_u - \alpha_d \geq -\alpha_{uu}$.  Hence we find that
\begin{align}
	\alpha_u - \alpha_d &\geq \left\lvert\alpha_{uu}\right\rvert \label{eqn:uuConstraint} \text{ .}
	\intertext{When the same analysis is performed with four cells, $\alpha_{uuu}$, $\alpha_{uu}$, $\alpha_u$ and $\alpha_d$, by varying $k \Delta x$ we find that equation~\eqref{eqn:uuConstraint} holds replacing $\left\lvert \alpha_{uu} \right\rvert$ with $\max(\left\lvert\alpha_{uu}\right\rvert, \left\lvert\alpha_{uuu}\right\rvert)$.  Hence, we generalise equation~\eqref{eqn:uuConstraint} to find the final stability constraint}
	\alpha_u - \alpha_d &\geq \max_{p\:\in\:P} |\alpha_p|
\end{align}
where the peripheral cells $P$ is the set of all stencil cells except for the upwind cell and downwind cell, and $\alpha_p$ is the weight for a given peripheral cell $p$.
We hypothesise that the three stability constraints (equations~\ref{eqn:vn:upwind}, \ref{eqn:vn:downwind} and \ref{eqn:uuConstraint}) are necessary but not sufficient for a transport scheme on arbitrary meshes.

\revtwo{The stability of the one-dimensional transport equation discretised in space and time could be analysed using existing techniques \citep{baldauf2008}, but we have only analysed the spatial stability of the cubicFit scheme.  Numerical experiments presented in section~\ref{sec:mountainAdvection} demonstrate that the cubicFit scheme is generally insensitive to the time-step, provided that it is below a stability limit.}

%% file: spherical.tex
\section*{Appendix B: Mesh geometry on a spherical Earth}

The cubicFit transport scheme is implemented using the OpenFOAM CFD library.  Unlike many atmospheric models that use spherical coordinates, OpenFOAM uses global, three-dimensional Cartesian coordinates with the $z$-axis pointing up through the North pole.  In order to perform the experiments on a spherical Earth presented in section~\ref{sec:deformationSphere}, it is necessary for velocity fields and mesh geometries to be expressed in these global Cartesian coordinates.

\subsection*{Velocity field specification}
The non-divergent velocity field in section~\ref{sec:deformationSphere} is specified as a streamfunction $\Psi(\lambda, \theta)$.  Instead of calculating velocity vectors, the flux $\mathbf{u}_f \cdot \mathbf{S}_f$ through a face $f$ is calculated directly from the streamfunction,
\begin{align}
	\mathbf{u}_f \cdot \mathbf{S}_f	= \sum_{e\:\in\:f} \mathbf{e} \cdot \mathbf{x}_e \Psi(e) \label{eqn:nondiv-spherical-flux}
\end{align}
where $e \in f$ denotes the edges $e$ of face $f$, $\mathbf{e}$ is the edge vector joining the two vertices of the edge, $\mathbf{x}_e$ is the position vector of the edge midpoint, and $\Psi(e)$ is the streamfunction evaluated at the same position.
Edge vectors are directed in a counter-clockwise orientation.

\subsection*{Spherical mesh construction}

Since OpenFOAM does not support two-dimensional spherical meshes, instead, we construct meshes that have a single layer of cells that are \SI{2000}{\meter} deep, having an inner radius $r_1 = R_e - \SI{1000}{\meter}$ and an outer radius $r_2 = R_e + \SI{1000}{\meter}$.
By default, OpenFOAM meshes comprise polyhedral cells with straight edges and flat faces.  This is problematic for spherical meshes because face areas and cell volumes are too small.
For tests on a spherical Earth, we override the default configuration and calculate our own face areas, cell volumes, face centres and cell centres that account for the mesh curvature.  Note that the new centres are no longer centroids, but they are consistent with the horizontal transport tests on a sphere presented in section~\ref{sec:deformationSphere}.

A face is classified as either a surface face or radial face.
A surface face has any number of vertices, all of equal radius.
A radial face has four vertices with two different radii, $r_1$ and $r_2$, and two different horizontal coordinates, $(\lambda_1, \theta_1)$ and $(\lambda_2, \theta_2)$.
A radial face centre is modified so that it has a radius $R_e$.  The latitudinal and longitudinal components of a radial face centre need no modification.
The face area $A_f$ for a radial face $f$ is the area of the annular sector,
\begin{align}
	A_f = \frac{d}{2} \left\lvert r_2^2 - r_1^2 \right\rvert
\end{align}
where $d$ is the great-circle distance between $(\lambda_1, \theta_1)$ and $(\lambda_2, \theta_2)$.

To calculate the centre of a surface face $f$, a new vertex is created that is positioned at the mean of the face vertices.  Note that this centre position, $\mathbf{\tilde{c}}_f$, is used in intermediate calculations and it is not the face centre position.
Next, the surface face is subdivided into spherical triangles that share this new vertex \citep{vanbrummelen2013}.
The face centre direction and radius are calculated separately.  The face centre direction $\mathbf{\hat{r}}$ is the mean of the spherical triangle centres weighted by their solid angle,
\begin{align}
	\mathbf{\hat{r}} = \frac
	{\sum_{t\:\in\:f}{\Omega_t \left(\mathbf{x}_{t,1} + \mathbf{x}_{t,2} + \mathbf{\tilde{c}}_f \right)}}
	{\left\lvert \sum_{t\:\in\:f}{\Omega_t \left(\mathbf{x}_{t,1} + \mathbf{x}_{t,2} + \mathbf{\tilde{c}}_f \right)} \right\rvert} \label{eqn:face-centre-dir}
\end{align}
where $t\in f$ denotes the spherical triangles $t$ of face $f$, $\Omega_t$ is spherical triangle's solid angle which is calculated using l'Huilier's theorem, $\mathbf{x}_{t,1}$ and $\mathbf{x}_{t,2}$ are the positions of the vertices shared by the face $f$ and spherical triangle $t$, and $\mathbf{\tilde{c}}_f$ is the position of the centre vertex shared by all spherical triangles of face $f$.
The face centre radius $r$ is the mean radius of the face vertices, again weighted by the solid angle of each spherical triangle,
\begin{align}
	r = \frac
	{\sum_{t\in f}{\Omega_t \left(\left\lvert \mathbf{x}_{t,1} \right\rvert + \left\lvert \mathbf{x}_{t,2} \right\rvert \right)/2}}
	{\Omega_f} \label{eqn:face-centre-mag}
\end{align}
where the solid angle $\Omega_f$ of face $f$ is the sum of the solid angles of the constituent spherical triangles,
\begin{align}
	\Omega_f = \sum_{t\in f}{\Omega_t} \text{ .}
\end{align}
We use equations~\eqref{eqn:face-centre-dir} and \eqref{eqn:face-centre-mag} to calculate the centre $\mathbf{c}_f$ of the face $f$,
\begin{align}
	\mathbf{c}_f = r\:\mathbf{\hat{r}} \label{eqn:face-centre}
\end{align}
The area vector $\mathbf{S}_f$ of the surface face $f$ is the sum of the spherical triangle areas \citep{vanbrummelen2013},
\begin{align}
	\mathbf{S}_f = r^2 \Omega_f \mathbf{\hat{r}} \text{ .}
\end{align}
Cell centres and cell volumes are corrected by considering faces that are not normal to the sphere such that
\begin{align}
	\frac{\left(\mathbf{S}_f \cdot \mathbf{c}_f\right)^2}{\left\lvert \mathbf{S}_f \right\rvert^2 \left\lvert \mathbf{c}_f \right\rvert^2} > 0 \label{eqn:surface-faces} \text{ .}
\end{align}
Let $\mathcal{F}$ be the set of faces satisfying equation~\eqref{eqn:surface-faces}.  Then, the cell volume $\mathcal{V}_c$ is
\begin{align}
	\mathcal{V}_c = \frac{1}{3} \sum_{f\:\in\:\mathcal{F}} \mathbf{S}_f \cdot \mathbf{c}_f
\end{align}
which can be thought of as the area $A$ integrated between $r_1$ and $r_2$ such that 
$\int_0^R{A(r)\:\mathrm{d}r} = \int_{r_1}^{r_2}{r^2 \Omega\:\mathrm{d}r} = \frac{1}{3} \Omega \left( r_2^3 - r_1^3 \right)$.
The cell centre is modified so that it has a radius $R_e$, which is consistent with radial faces.

Edges can be classified in a similar manner to faces where surface edges are tangent to the sphere and radial faces are normal to the sphere.  The edge midpoints $\mathbf{x}_e$ are used to calculate the face flux for non-divergent velocity fields (equation~\ref{eqn:nondiv-spherical-flux}).
For transport tests, corrections to edge midpoints are unnecessary.  Due to the choice of $r_1$ and $r_2$ during mesh construction, the midpoint of a radial edge is at a radial distance of $R_e$ which is necessary for the correct calculation of non-divergent velocity fields.
The position of surface edge midpoints is unimportant because these edges do not contribute to the face flux since $\mathbf{e} \cdot \mathbf{x}_e = 0$.
Edge lengths are the straight-line distance between the two vertices and not the great-circle distance.  Again, the edge lengths are not corrected because it makes no difference to the face flux calculation.

%% file: slant-cell-advection.bbl
\begin{thebibliography}{10}
\expandafter\ifx\csname url\endcsname\relax
  \def\url#1{\texttt{#1}}\fi
\expandafter\ifx\csname urlprefix\endcsname\relax\def\urlprefix{URL }\fi
\expandafter\ifx\csname href\endcsname\relax
  \def\href#1#2{#2} \def\path#1{#1}\fi

\bibitem{walko-avissar2008b}
R.~L. Walko, R.~Avissar, The {Ocean-Land-Atmosphere Model (OLAM)}. {Part II:
  Formulation} and tests of the nonhydrostatic dynamic core, Mon.\ Wea.\ Rev.
  136~(11) (2008) 4045--4062.
\newblock \href {http://dx.doi.org/10.1175/2008MWR2523.1}
  {\path{doi:10.1175/2008MWR2523.1}}.

\bibitem{schaer2002}
C.~Sch{\"a}r, D.~Leuenberger, O.~Fuhrer, D.~L{\"u}thi, C.~Girard, A new
  terrain-following vertical coordinate formulation for atmospheric prediction
  models, Mon.\ Wea.\ Rev. 130~(10) (2002) 2459--2480.
\newblock \href
  {http://dx.doi.org/10.1175/1520-0493(2002)130<2459:ANTFVC>2.0.CO;2}
  {\path{doi:10.1175/1520-0493(2002)130<2459:ANTFVC>2.0.CO;2}}.

\bibitem{hoinka-zaengl2004}
K.~P. Hoinka, G.~Z{\"a}ngl, The influence of the vertical coordinate on
  simulations of a pv streamer crossing the alps, Mon.\ Wea.\ Rev. 132~(7)
  (2004) 1860--1867.
\newblock \href
  {http://dx.doi.org/10.1175/1520-0493(2004)132<1860:TIOTVC>2.0.CO;2}
  {\path{doi:10.1175/1520-0493(2004)132<1860:TIOTVC>2.0.CO;2}}.

\bibitem{webster2003}
S.~Webster, A.~Brown, D.~Cameron, C.~Jones, Improvements to the representation
  of orography in the {Met Office Unified Model}, Quart.\ J.\ Roy.\ Meteor.\
  Soc. 129~(591) (2003) 1989--2010.
\newblock \href {http://dx.doi.org/10.1256/qj.02.133}
  {\path{doi:10.1256/qj.02.133}}.

\bibitem{galchen-somerville1975}
T.~Gal-Chen, R.~C. Somerville, On the use of a coordinate transformation for
  the solution of the {Navier-Stokes} equations, J.\ Comp.\ Phys. 17~(2) (1975)
  209--228.
\newblock \href {http://dx.doi.org/10.1016/0021-9991(75)90037-6}
  {\path{doi:10.1016/0021-9991(75)90037-6}}.

\bibitem{klemp2011}
J.~B. Klemp, A terrain-following coordinate with smoothed coordinate surfaces,
  Mon.\ Wea.\ Rev. 139~(7) (2011) 2163--2169.
\newblock \href {http://dx.doi.org/10.1175/MWR-D-10-05046.1}
  {\path{doi:10.1175/MWR-D-10-05046.1}}.

\bibitem{eckermann2014}
S.~D. Eckermann, J.~P. McCormack, J.~Ma, T.~F. Hogan, K.~A. Zawdie,
  Stratospheric analysis and forecast errors using hybrid and sigma
  coordinates, Mon.\ Wea.\ Rev. 142~(1) (2014) 476--485.
\newblock \href {http://dx.doi.org/10.1175/MWR-D-13-00203.1}
  {\path{doi:10.1175/MWR-D-13-00203.1}}.

\bibitem{simmons-burridge1981}
A.~J. Simmons, D.~M. Burridge, An energy and angular-momentum conserving
  vertical finite-difference scheme and hybrid vertical coordinates, Mon.\
  Wea.\ Rev. 109~(4) (1981) 758--766.
\newblock \href
  {http://dx.doi.org/10.1175/1520-0493(1981)109<0758:AEAAMC>2.0.CO;2}
  {\path{doi:10.1175/1520-0493(1981)109<0758:AEAAMC>2.0.CO;2}}.

\bibitem{leuenberger2010}
D.~Leuenberger, M.~Koller, O.~Fuhrer, C.~Sch{\"a}r, A generalization of the
  {SLEVE} vertical coordinate, Mon.\ Wea.\ Rev. 138~(9) (2010) 3683--3689.
\newblock \href {http://dx.doi.org/10.1175/2010MWR3307.1}
  {\path{doi:10.1175/2010MWR3307.1}}.

\bibitem{yamazaki2016}
H.~Yamazaki, T.~Satomura, N.~Nikiforakis, Three-dimensional cut-cell modelling
  for high-resolution atmospheric simulations, Quart.\ J.\ Roy.\ Meteor.\ Soc.
  142~(696) (2016) 1335--1350.
\newblock \href {http://dx.doi.org/10.1002/qj.2736}
  {\path{doi:10.1002/qj.2736}}.

\bibitem{steppeler2002}
J.~Steppeler, H.-W. Bitzer, M.~Minotte, L.~Bonaventura, Nonhydrostatic
  atmospheric modeling using a $z$-coordinate representation, Mon.\ Wea.\ Rev.
  130~(8) (2002) 2143--2149.
\newblock \href
  {http://dx.doi.org/10.1175/1520-0493(2002)130<2143:NAMUAZ>2.0.CO;2}
  {\path{doi:10.1175/1520-0493(2002)130<2143:NAMUAZ>2.0.CO;2}}.

\bibitem{jebens2011}
S.~Jebens, O.~Knoth, R.~Weiner, Partially implicit peer methods for the
  compressible {Euler} equations, J.\ Comp.\ Phys. 230~(12) (2011) 4955--4974.
\newblock \href {http://dx.doi.org/10.1016/j.jcp.2011.03.015}
  {\path{doi:10.1016/j.jcp.2011.03.015}}.

\bibitem{shaw-weller2016}
J.~Shaw, H.~Weller, Comparison of terrain following and cut cell grids using a
  non-hydrostatic model, Mon.\ Wea.\ Rev. 144~(6) (2016) 2085--2099.
\newblock \href {http://dx.doi.org/10.1175/MWR-D-15-0226.1}
  {\path{doi:10.1175/MWR-D-15-0226.1}}.

\bibitem{leonard1993}
B.~Leonard, M.~MacVean, A.~Lock, Positivity-preserving numerical schemes for
  multidimensional advection, Tech. Rep. 106055, NASA (1993).

\bibitem{kent2014}
J.~Kent, P.~A. Ullrich, C.~Jablonowski, Dynamical core model intercomparison
  project: Tracer transport test cases, Quart.\ J.\ Roy.\ Meteor.\ Soc.
  140~(681) (2014) 1279--1293.
\newblock \href {http://dx.doi.org/10.1002/qj.2208}
  {\path{doi:10.1002/qj.2208}}.

\bibitem{leonard1996}
B.~Leonard, A.~Lock, M.~MacVean, Conservative explicit unrestricted-time-step
  multidimensional constancy-preserving advection schemes, Mon.\ Wea.\ Rev.
  124~(11) (1996) 2588--2606.
\newblock \href
  {http://dx.doi.org/10.1175/1520-0493(1996)124<2588:CEUTSM>2.0.CO;2}
  {\path{doi:10.1175/1520-0493(1996)124<2588:CEUTSM>2.0.CO;2}}.

\bibitem{weller2017}
Y.~Chen, H.~Weller, S.~Pring, J.~Shaw, Dimension splitting errors and a long
  time-step multi-dimensional scheme for atmospheric transport, Quart.\ J.\
  Roy.\ Meteor.\ Soc.In revision; preprint,
  \url{http://arxiv.org/abs/1701.06907}.

\bibitem{lin-rood1996}
S.-J. Lin, R.~B. Rood, Multidimensional flux-form {semi-Lagrangian} transport
  schemes, Mon.\ Wea.\ Rev. 124~(9) (1996) 2046--2070.
\newblock \href
  {http://dx.doi.org/10.1175/1520-0493(1996)124<2046:MFFSLT>2.0.CO;2}
  {\path{doi:10.1175/1520-0493(1996)124<2046:MFFSLT>2.0.CO;2}}.

\bibitem{guo2014}
W.~Guo, R.~D. Nair, J.-M. Qiu, A conservative {semi-Lagrangian} discontinuous
  {Galerkin} scheme on the cubed sphere, Mon.\ Wea.\ Rev. 142~(1) (2014)
  457--475.
\newblock \href {http://dx.doi.org/10.1175/MWR-D-13-00048.1}
  {\path{doi:10.1175/MWR-D-13-00048.1}}.

\bibitem{putman-lin2007}
W.~M. Putman, S.-J. Lin, Finite-volume transport on various cubed-sphere grids,
  J.\ Comp.\ Phys. 227~(1) (2007) 55--78.
\newblock \href {http://dx.doi.org/10.1016/j.jcp.2007.07.022}
  {\path{doi:10.1016/j.jcp.2007.07.022}}.

\bibitem{iske-kaeser2004}
A.~Iske, M.~K{\"a}ser, Conservative {semi-Lagrangian} advection on adaptive
  unstructured meshes, Numer.\ Methods.\ Partial\ Differ.\ Equ. 20~(3) (2004)
  388--411.
\newblock \href {http://dx.doi.org/10.1002/num.10100}
  {\path{doi:10.1002/num.10100}}.

\bibitem{lauritzen2010}
P.~H. Lauritzen, R.~D. Nair, P.~A. Ullrich, A conservative semi-lagrangian
  multi-tracer transport scheme ({CSLAM}) on the cubed-sphere grid, J.\ Comp.\
  Phys. 229~(5) (2010) 1401--1424.
\newblock \href {http://dx.doi.org/10.1016/j.jcp.2009.10.036}
  {\path{doi:10.1016/j.jcp.2009.10.036}}.

\bibitem{lauritzen2011book}
P.~H. Lauritzen, C.~Jablonowski, M.~A. Taylor, R.~D. Nair, Numerical techniques
  for global atmospheric models, Vol.~80, Springer Science \& Business Media,
  2011.
\newblock \href {http://dx.doi.org/10.1007/978-3-642-11640-7}
  {\path{doi:10.1007/978-3-642-11640-7}}.

\bibitem{zerroukat-allen2012}
M.~Zerroukat, T.~Allen, A three-dimensional monotone and conservative
  {semi-Lagrangian} scheme ({SLICE-3D}) for transport problems, Quart.\ J.\
  Roy.\ Meteor.\ Soc. 138~(667) (2012) 1640--1651.
\newblock \href {http://dx.doi.org/10.1002/qj.1902}
  {\path{doi:10.1002/qj.1902}}.

\bibitem{lashley2002}
R.~K. Lashley, Automatic generation of accurate advection schemes on
  unstructured grids and their application to meteorological problems, Ph.D.
  thesis, University of Reading (2002).

\bibitem{skamarock-menchaca2010}
W.~C. Skamarock, M.~Menchaca, Conservative transport schemes for spherical
  geodesic grids: High-order reconstructions for forward-in-time schemes, Mon.\
  Wea.\ Rev. 138~(12) (2010) 4497--4508.
\newblock \href {http://dx.doi.org/10.1175/2010MWR3390.1}
  {\path{doi:10.1175/2010MWR3390.1}}.

\bibitem{lauritzen2011}
P.~H. Lauritzen, C.~Erath, R.~Mittal, On simplifying ‘incremental
  remap’--based transport schemes, J.\ Comp.\ Phys. 230~(22) (2011)
  7957--7963.
\newblock \href {http://dx.doi.org/10.1016/j.jcp.2011.06.030}
  {\path{doi:10.1016/j.jcp.2011.06.030}}.

\bibitem{thuburn2014}
J.~Thuburn, C.~Cotter, T.~Dubos, A mimetic, semi-implicit, forward-in-time,
  finite volume shallow water model: comparison of hexagonal-icosahedral and
  cubed-sphere grids, Geosci.\ Model\ Dev. 7~(3) (2014) 909--929.
\newblock \href {http://dx.doi.org/10.5194/gmd-7-909-2014}
  {\path{doi:10.5194/gmd-7-909-2014}}.

\bibitem{weller2009}
H.~Weller, H.~G. Weller, A.~Fournier, {Voronoi}, {Delaunay}, and
  block-structured mesh refinement for solution of the shallow-water equations
  on the sphere, Mon.\ Wea.\ Rev. 137~(12) (2009) 4208--4224.
\newblock \href {http://dx.doi.org/10.1175/2009MWR2917.1}
  {\path{doi:10.1175/2009MWR2917.1}}.

\bibitem{skamarock-gassmann2011}
W.~C. Skamarock, A.~Gassmann, Conservative transport schemes for spherical
  geodesic grids: {High-order} flux operators for {ODE}-based time integration,
  Mon.\ Wea.\ Rev. 139~(9) (2011) 2962--2975.
\newblock \href {http://dx.doi.org/10.1175/MWR-D-10-05056.1}
  {\path{doi:10.1175/MWR-D-10-05056.1}}.

\bibitem{steppeler-klemp2017}
J.~Steppeler, J.~Klemp, Advection on cut-cell grids for an idealized mountain
  of constant slope, Mon.\ Wea.\ Rev.\href
  {http://dx.doi.org/10.1175/MWR-D-16-0308.1}
  {\path{doi:10.1175/MWR-D-16-0308.1}}.

\bibitem{lauritzen2014}
P.~Lauritzen, P.~Ullrich, C.~Jablonowski, P.~Bosler, D.~Calhoun, A.~Conley,
  T.~Enomoto, L.~Dong, S.~Dubey, O.~Guba, et~al., A standard test case suite
  for two-dimensional linear transport on the sphere: results from a collection
  of state-of-the-art schemes, Geosci.\ Model\ Dev. 7~(1) (2014) 105--145.
\newblock \href {http://dx.doi.org/10.5194/gmd-7-105-2014}
  {\path{doi:10.5194/gmd-7-105-2014}}.

\bibitem{ii-xiao2009}
S.~Ii, F.~Xiao, High order multi-moment constrained finite volume method. {Part
  I}: Basic formulation, J.\ Comp.\ Phys. 228~(10) (2009) 3669--3707.
\newblock \href {http://dx.doi.org/10.1016/j.jcp.2009.02.009}
  {\path{doi:10.1016/j.jcp.2009.02.009}}.

\bibitem{li2013}
X.~Li, C.~Chen, X.~Shen, F.~Xiao, A multimoment constrained finite-volume model
  for nonhydrostatic atmospheric dynamics, Mon.\ Wea.\ Rev. 141~(4) (2013)
  1216--1240.
\newblock \href {http://dx.doi.org/10.1175/MWR-D-12-00144.1}
  {\path{doi:10.1175/MWR-D-12-00144.1}}.

\bibitem{xie-xiao2016}
B.~Xie, F.~Xiao, A multi-moment constrained finite volume method on arbitrary
  unstructured grids for incompressible flows, J.\ Comp.\ Phys. 327 (2016)
  747--778.

\bibitem{smolarkiewicz-szmelter2011}
P.~K. Smolarkiewicz, J.~Szmelter, A nonhydrostatic unstructured-mesh soundproof
  model for simulation of internal gravity waves, Acta\ Geophys. 59~(6) (2011)
  1109--1134.
\newblock \href {http://dx.doi.org/10.2478/s11600-011-0043-z}
  {\path{doi:10.2478/s11600-011-0043-z}}.

\bibitem{smolarkiewicz-szmelter2005}
P.~K. Smolarkiewicz, J.~Szmelter, {MPDATA}: An edge-based unstructured-grid
  formulation, J.\ Comp.\ Phys. 206~(2) (2005) 624--649.
\newblock \href {http://dx.doi.org/10.1016/j.jcp.2004.12.021}
  {\path{doi:10.1016/j.jcp.2004.12.021}}.

\bibitem{kuehnlein-smolarkiewicz2017}
C.~K{\"u}hnlein, P.~K. Smolarkiewicz, An unstructured-mesh finite-volume
  {MPDATA} for compressible atmospheric dynamics, J.\ Comp.\ Phys. 334 (2017)
  16--30.
\newblock \href {http://dx.doi.org/10.1016/j.jcp.2016.12.054}
  {\path{doi:10.1016/j.jcp.2016.12.054}}.

\bibitem{weller-shahrokhi2014}
H.~Weller, A.~Shahrokhi, Curl free pressure gradients over orography in a
  solution of the fully compressible {Euler} equations with implicit treatment
  of acoustic and gravity waves, Mon.\ Wea.\ Rev. 142~(12) (2014) 4439--4457.
\newblock \href {http://dx.doi.org/10.1175/MWR-D-14-00054.1}
  {\path{doi:10.1175/MWR-D-14-00054.1}}.

\bibitem{cuetofelgueroso2006}
L.~Cueto-Felgueroso, I.~Colominas, J.~Fe, F.~Navarrina, M.~Casteleiro,
  High-order finite volume schemes on unstructured grids using moving
  least-squares reconstruction. {Application} to shallow water dynamics, Int.\
  J.\ Numer.\ Meth.\ Engng 65~(3) (2006) 295--331.
\newblock \href {http://dx.doi.org/10.1002/nme.1442}
  {\path{doi:10.1002/nme.1442}}.

\bibitem{cuetofelgueroso2007}
L.~Cueto-Felgueroso, I.~Colominas, X.~Nogueira, F.~Navarrina, M.~Casteleiro,
  Finite volume solvers and moving least-squares approximations for the
  compressible {Navier--Stokes} equations on unstructured grids, Comput.\
  Methods\ Appl.\ Mech.\ Engrg 196~(45) (2007) 4712--4736.
\newblock \href {http://dx.doi.org/10.1016/j.cma.2007.06.003}
  {\path{doi:10.1016/j.cma.2007.06.003}}.

\bibitem{white2017}
R.~P. Laurent~White, D.~Trenev, Flow simulation in heterogeneous porous media
  with the moving least-squares method, SIAM J.\ Sci.\ Comput. 39~(2) (2017)
  B323--B351.
\newblock \href {http://dx.doi.org/10.1137/16M1070840}
  {\path{doi:10.1137/16M1070840}}.

\bibitem{lauritzen2012}
P.~H. Lauritzen, W.~C. Skamarock, M.~Prather, M.~Taylor, A standard test case
  suite for two-dimensional linear transport on the sphere, Geosci.\ Model\
  Dev. 5~(3) (2012) 887--901.
\newblock \href {http://dx.doi.org/10.5194/gmd-5-887-2012}
  {\path{doi:10.5194/gmd-5-887-2012}}.

\bibitem{durran2013}
D.~R. Durran, Numerical methods for wave equations in geophysical fluid
  dynamics, Vol.~32, Springer Science \& Business Media, 2013.
\newblock \href {http://dx.doi.org/10.1007/978-1-4419-6412-0}
  {\path{doi:10.1007/978-1-4419-6412-0}}.

\bibitem{openfoam}
{CFD Direct}, {OpenFOAM} user guide: Numerical schemes,
  \url{http://cfd.direct/openfoam/user-guide/fvschemes/} (2016).

\bibitem{atmosfoam}
H.~Weller, J.~Shaw, {AtmosFOAM}: {OpenFOAM} applications and libraries for
  atmospheric modelling (2017).
\newblock \href {http://dx.doi.org/10.5281/zenodo.247327}
  {\path{doi:10.5281/zenodo.247327}}.

\bibitem{asam_grid}
{Leibniz Institute for Tropospheric Research}, J.~Shaw, {ASAM} cut cell grid
  generator (2017).
\newblock \href {http://dx.doi.org/10.5281/zenodo.242374}
  {\path{doi:10.5281/zenodo.242374}}.

\bibitem{gmv2openfoam}
J.~Shaw, {GMV-to-OpenFOAM} file format converter (2017).
\newblock \href {http://dx.doi.org/10.5281/zenodo.242387}
  {\path{doi:10.5281/zenodo.242387}}.

\bibitem{geodesic-mesh}
J.~Thuburn, C.~J. Cotter, T.~Dubos, H.~Weller, J.~Shaw, {Hexagonal/triangular
  geodesic mesh generator} (2017).
\newblock \href {http://dx.doi.org/10.5281/zenodo.245327}
  {\path{doi:10.5281/zenodo.245327}}.

\bibitem{atmostests}
J.~Shaw, {AtmosTests}: Idealised atmospheric test suite (2017).
\newblock \href {http://dx.doi.org/10.5281/zenodo.247480}
  {\path{doi:10.5281/zenodo.247480}}.

\bibitem{atmostests-data}
J.~Shaw, Test result data for idealised numerical atmospheric transport (2017).
\newblock \href {http://dx.doi.org/10.5281/zenodo.259328}
  {\path{doi:10.5281/zenodo.259328}}.

\bibitem{good2014}
B.~Good, A.~Gadian, S.-J. Lock, A.~Ross, Performance of the cut-cell method of
  representing orography in idealized simulations, Atmos. Sci. Lett. 15~(1)
  (2014) 44--49.
\newblock \href {http://dx.doi.org/10.1002/asl2.465}
  {\path{doi:10.1002/asl2.465}}.

\bibitem{jaehn2015}
M.~J\"ahn, O.~Knoth, M.~K\"onig, U.~Vogelsberg, {ASAM v2.7}: a compressible
  atmospheric model with a {Cartesian} cut cell approach, Geosci.\ Model\ Dev.
  8~(2) (2015) 317--340.
\newblock \href {http://dx.doi.org/10.5194/gmd-8-317-2015}
  {\path{doi:10.5194/gmd-8-317-2015}}.

\bibitem{asam2010}
U.~Vogelsberg, {ASAM} wiki: Grid generator,
  \url{http://asamwiki.tropos.de/index.php?title=Grid_Generator} (2010).

\bibitem{baldauf2008}
M.~Baldauf, Stability analysis for linear discretisations of the advection
  equation with {Runge--Kutta} time integration, J.\ Comp.\ Phys. 227~(13)
  (2008) 6638--6659.
\newblock \href {http://dx.doi.org/10.1016/j.jcp.2008.03.025}
  {\path{doi:10.1016/j.jcp.2008.03.025}}.

\bibitem{heikes-randall1995a}
R.~Heikes, D.~A. Randall, Numerical integration of the shallow-water equations
  on a twisted icosahedral grid. {Part I}: Basic design and results of tests,
  Mon.\ Wea.\ Rev. 123~(6) (1995) 1862--1880.
\newblock \href
  {http://dx.doi.org/10.1175/1520-0493(1995)123<1862:NIOTSW>2.0.CO;2}
  {\path{doi:10.1175/1520-0493(1995)123<1862:NIOTSW>2.0.CO;2}}.

\bibitem{heikes-randall1995b}
R.~Heikes, D.~A. Randall, Numerical integration of the shallow-water equations
  on a twisted icosahedral grid. {Part II}: A detailed description of the grid
  and an analysis of numerical accuracy, Mon.\ Wea.\ Rev. 123~(6) (1995)
  1881--1887.
\newblock \href
  {http://dx.doi.org/10.1175/1520-0493(1995)123<1881:NIOTSW>2.0.CO;2}
  {\path{doi:10.1175/1520-0493(1995)123<1881:NIOTSW>2.0.CO;2}}.

\bibitem{staniforth-thuburn2012}
A.~Staniforth, J.~Thuburn, Horizontal grids for global weather and climate
  prediction models: a review, Quart.\ J.\ Roy.\ Meteor.\ Soc. 138~(662) (2012)
  1--26.
\newblock \href {http://dx.doi.org/10.1002/qj.958} {\path{doi:10.1002/qj.958}}.

\bibitem{white-dongarra2011}
J.~White, J.~J. Dongarra, High-performance high-resolution {semi-Lagrangian}
  tracer transport on a sphere, J.\ Comp.\ Phys. 230~(17) (2011) 6778--6799.
\newblock \href {http://dx.doi.org/10.1016/j.jcp.2011.05.008}
  {\path{doi:10.1016/j.jcp.2011.05.008}}.

\bibitem{miura2007}
H.~Miura, An upwind-biased conservative advection scheme for spherical
  hexagonal-pentagonal grids, Mon.\ Wea.\ Rev. 135~(12) (2007) 4038--4044.
\newblock \href {http://dx.doi.org/10.1175/2007MWR2101.1}
  {\path{doi:10.1175/2007MWR2101.1}}.

\bibitem{vanbrummelen2013}
G.~Van~Brummelen, Heavenly mathematics: The forgotten art of spherical
  trigonometry, Princeton University Press, 2013.

\end{thebibliography}
